\documentclass[11pt,a4paper]{article}

\usepackage[utf8]{inputenc}  
\usepackage{amsmath, amssymb}  
\usepackage{graphicx}  
\usepackage{hyperref}  
\usepackage{geometry}  
\usepackage{color}  
\usepackage{xcolor}          
\usepackage{cite}  
\usepackage{authblk}  
\usepackage{booktabs}
\usepackage{tabularx}

\usepackage{lipsum}
\usepackage{amsfonts}
\usepackage{subfigure}
\usepackage{float}    
\usepackage{comment}
\usepackage{epstopdf}
\usepackage{algorithmic}
\usepackage{bbm}
\usepackage{textgreek}	

\title{ Well-posedness and finite-time breakdown for Euler-alignment systems with diffusive chemotaxis}

\author{
Marta Menci\thanks{Department of Engineering, Università Campus Bio-Medico di Roma, Rome, Italy.
E-mail: m.menci@unicampus.it},$\;\;$ Marco Papi\thanks{Department of Engineering, Università Campus Bio-Medico di Roma, Rome, Italy.
E-mail: m.papi@unicampus.it}
}
\date{}

\newtheorem{theorem}{Theorem}
\newtheorem{proposition}{Proposition}
\newtheorem{lemma}{Lemma}
\newtheorem{remark}{Remark}
\newtheorem{proof}{Proof}
\newtheorem{Definition}{Definition}

\begin{document}

\maketitle

\begin{abstract}\noindent
We investigate a recent class of pressureless Euler systems arising in biological contexts, in which the velocity evolution is driven, besides the nonlocal alignment interaction, 
by the gradient of a concentration self-produced by the population. 
We first establish local existence and uniqueness of strong solutions, in a class requiring no spatial derivative of the initial population density, together with a continuation criterion. We then derive sufficient conditions on the initial
data for finite-time breakdown in the multidimensional setting, with an explicit upper bound for the maximal time of existence.
Numerical simulations complement the analysis and illustrate the transition between regular and singular regimes, and the role of the chemotactic strength in the onset of blow-up.

\end{abstract}


\medskip


\medskip

\noindent\textbf{Corresponding author:} Marta Menci (\texttt{m.menci@unicampus.it})

\section{Introduction}
\label{sec:new_intro}

In this paper we investigate finite-time breakdown and  well-posedness for a generalized version of pressureless
Euler alignment systems coupling non-local velocity alignment and chemotaxis.
Denoting with $\rho$ a nonnegative density, $u$ the velocity field, and $f$ the concentration of the chemical signal influencing the dynamics, the system at a glance reads 
\begin{equation}
\label{eq:model_intro}
\begin{cases}
\displaystyle
\partial_t\rho+\nabla\cdot(\rho u)=0,
\\[6pt]
\displaystyle
\partial_t u+(u\cdot\nabla)u
=\int_{\mathbb R^n}\phi(|x-y|)\big(u(y,t)-u(x,t)\big)\rho(y,t)\,dy
+\gamma\,\nabla f,
\\[8pt]
\displaystyle
\partial_t f-\operatorname{tr}\!\big(\sigma\sigma^\top\,\nabla^2 f\big)+b\cdot\nabla f
=\psi*\rho,
\end{cases}
\qquad x\in\mathbb R^n,\ \ t\in(0,T),
\end{equation}
subject to the initial data $\rho(\cdot,0)=\rho_0$, $u(\cdot,0)=u_0$, and $f(\cdot,0)=f_0$.

The system shares common features of pressureless Euler alignment systems in \cite{carrillo2017review, albi2018pressureless, tadmor2014critical, choi2017emergent}, where the equation for the velocity is driven by a
nonlocal alignment term. Here, the velocity at $x$ is relaxed towards the velocities at neighbouring points $y$, with weights determined by the density $\rho$ and by a given \emph{communication kernel} $\phi$. 
A key aspect of \eqref{eq:model_intro} relies on a chemotactic drift $\gamma\,\nabla f$, steering the population dynamics towards regions of higher chemoattractant concentration if $\gamma>0$.
The evolution of the chemical concentration $f$ is ruled by a parabolic equation, which here generalizes the structure of the equation considered in hybrid-coupled models of the literature \cite{di2015hybrid, costanzo2020hybrid, menci2023microscopic}.
In particular, the equation models collective cell dynamics scenarios where the chemoattractant is produced by the cell population, whose density is here denoted with $\rho$, through the nonlocal source $\psi*\rho$, with $\psi$ denoting a localized smooth kernel referred as the \textit{chemoattractant kernel}.

Although the system is here motivated by biological applications, its mathematical structure admits a broader interpretation. For instance, in social and opinion dynamics settings, the chemical signal can be seen as a distributed information field generated by the population, which propagates and feeds back to the dynamics through its spatial gradient.  
The coupling with $f$ in the velocity equation introduces an additional mechanism, in which agents respond to an endogenous information environment created by the past evolution of the population. 

Throughout the rest of the paper we refer to \eqref{eq:model_intro}
as the \emph{Euler-alignment system with chemotaxis}. 
The novelty of the model under investigation mainly relies on the coupling between $(\rho, u)$ and $f$, i.e. the interplay between alignment and chemotaxis, which acts
as a history-dependent feedback in the velocity dynamics.
The results are obtained for the multidimensional space domain $\mathbb{R}^n$, with $n\ge1$, generalizing the structure of the parabolic equation with respect to the literature. 
The model in \eqref{eq:model_intro} arises as the macroscopic counterpart of a class of hybrid discrete-continuum models, originally conceived for modelling collective dynamics of cells under mechanical interactions and chemotaxis \cite{di2015hybrid, di2018discrete}.

\paragraph{From hybrid particle-continuum models to a hydrodynamic description}
\label{subsec:derivation}

In the hybrid description, a finite number $N$ of agents (cells) is described by positions $x_i(t)$ and velocities $v_i(t)$, whose dynamics are coupled to a parabolic diffusion equations, with constant coefficients:
\begin{equation}
\label{eq:particle_intro}
\begin{cases}
\displaystyle
\dot x_i=v_i,
\\[4pt]
\displaystyle
\dot v_i=\frac1N\sum_{j=1}^N\phi(|x_i-x_j|)\,(v_j-v_i)+\gamma\,\nabla f(x_i,t),
\\[6pt]
\displaystyle
\partial_t f=D \Delta f + g(x,X(t)),
\end{cases}
\qquad i=1,\dots,N.
\end{equation}
where $D, \gamma \in \mathbb{R}^{+}$.
Here the agents interact through the same alignment kernel $\phi$ and respond to the gradient of
a self-produced chemoattractant, the production being modelled by function $g$, that depends on the positions of the cells through $X:=(x_1,...,x_N)$.

Hybrid systems of this type have been analytically studied in 
\cite{menci2019global, menci2022existence, menci2023coupled, menci2023existence, costanzo2020hybrid}, establishing a
rigorous analytical framework for this kind of models. In this
setting, local and global well-posedness results were established under suitable regularity
assumptions on the source term \cite{menci2019global}, and subsequently extended to more general
settings, including discontinuous production terms \cite{menci2022existence, menci2023coupled}. 
Passing formally to the mean-field regime $N\to\infty$,
the alignment interaction in \eqref{eq:particle_intro} converges to the nonlocal integral term
$
\int_{\mathbb R^n}\phi(|x-y|)\big(u(y,t)-u(x,t)\big)\rho(y,t)\,dy.
$
The rigorous derivation of the macroscopic counterpart for the class of hybrid-coupled models has been performed in \cite{natalini2021mean}, and numerical studies for the models at the different scales of description have been performed in \cite{menci2023microscopic, menci2026numericalMMS}.

\paragraph{Critical thresholds and well-posedeness results without chemotaxis}
\label{subsec:no_chemo}

When the chemotactic coupling is switched off, i.e. $\gamma=0$ and removing the third equation in \eqref{eq:model_intro},  the pair $(\rho,u)$ solves the pressureless \emph{Euler-alignment system}
\begin{equation}
\label{eq:EAS_intro}
\begin{cases}
\partial_t\rho+\nabla\cdot(\rho u)=0,\\[4pt]
\partial_t u+(u\cdot\nabla)u=\displaystyle\int_{\mathbb R^n}\phi(|x-y|)\big(u(y,t)-u(x,t)\big)\rho(y,t)\,dy.
\end{cases}
\end{equation}

which describes, at the hydrodynamic level, collective motions driven by nonlocal velocity alignment. The literature concerning the rigorous derivation from the particle level through kinetic one is currently rich, see \cite{carrillo2021mean, golse2016dynamics, carrillo2026lagrangian, choi2025alignment, jabin2014review, carrillo2017review, choi2017emergent} and references therein for seminal and recent papers of the literature. 
Also in the case of topological interactions, thus where interactions do not depend on metric distance, the path from microscopic to macroscopic pressureless Euler-type system has been investigated, both analytically and numerically \cite{Haskovec2013FlockingDA, DP19, menci2026numerical}. 

From the analytical viewpoint, the system has been investigated in different directions.
In particular, the system has been investigated in the framework of \textit{critical thresholds}, establishing conditions leading to global solutions and asymptotic flocking, versus finite-time blow-up.
Starting from the work in \cite{tadmor2014critical}, several results follow in the one-dimensional and two-dimensional case, see for instance \cite{carrillo2016critical,carrillo2017review, global2D_Tadmor2017,albi2018pressureless}. 
Moreover, the study of critical threshold phenomena have been extended to other Euler-type systems, including the pressured Euler systems and general Euler-Poisson systems \cite{tadmor2003critical, luan2026euler, wei2012critical, bhatnagar2020critical, tadmor2008global}. Re

\paragraph{Paper contribution and structure}

Our contribution is threefold. In Section \ref{sec:wellpos}, we establish local well-posedness and a continuation
criterion for \eqref{eq:model_intro} under weak assumptions on the initial data, the kernels and the coefficient of the parabolic equations.
Then, a continuation criterion shows that finite-time breakdown may occur only through the blow-up of $\nabla u$; in particular, the
chemotactic drift does not destroy the a priori bound on $\|u\|_{L^\infty}$, which we obtain
explicitly in terms of the data.
In Section \ref{sec:multidimensional_thresholds}, working directly on the velocity-gradient matrix along the characteristics, we derive a
sufficient condition on the initial data for finite-time breakdown in the fully multidimensional
setting, together with an explicit upper bound for the maximal time of existence. 
In one space dimension we obtain a sufficient condition, verifiable on the initial velocity
slope, ensuring that the solution can be continued beyond a prescribed time, and we compare the
behavior with the well-investigated chemotaxis-free dynamics ($\gamma=0$). Moreover, in the case of constant parabolic coefficients the conditions become
fully explicit. In this setting, we validate numerically the results obtained, with a detailed study of the sensitivity of detected blow-up with respect to the discretization parameters and of the role of chemotaxis (Section \ref{sec:numerical_experiments}).
As stated in the recent paper \cite{choi2026finite}, yet for Euler-alignment systems quantitative criteria for finite-time blow-up are still less developed with respect to global results. This is even more evident for Euler-alignment system with chemotaxis, for which the present work represents a first step towards rigorous results and critical threshold framework.

\section{Local well-posedness and continuation of solutions } 
\label{sec:wellpos}

For the sake of clarity we first introduce the notation that will be used in the following.
Throughout the paper, for any $v,w\in\mathbb R^n$, $v\otimes w$ denotes the $n\times n$-matrix with entries $(v\otimes w)_{ij}=v_i w_j$.
Let $f:\mathbb{R}^n \to \mathbb{R}$ and a multi-index $\alpha=(\alpha_1,\dots,\alpha_n)\in\mathbb{N}^n$, we denote by
$
D_x^\alpha f = \partial_{x_1}^{\alpha_1}\cdots \partial_{x_n}^{\alpha_n} f
$
the corresponding partial derivative of order $|\alpha|=\alpha_1+\cdots+\alpha_n$, with respect to the spatial variable $x$, and by $\nabla f$ the spatial gradient of the scalar function $f$, $\nabla^2 f$ the Hessian matrix of $f$ with respect to the spatial variable $x$.

For the vector field $u=(u_1,\ldots,u_n)^\top$ $:\mathbb{R}^n \to \mathbb{R}^n$, we denote by $\nabla u$ its spatial gradient, defined as the Jacobian matrix
$
\nabla u(x)
:=
\left( \partial_{x_j} u_i(x) \right)_{1\le i,j\le n}
\in \mathbb{R}^{n\times n}.
$
In particular, the $(i,j)$-th entry of $\nabla u$ is given by $\partial_{x_j} u_i(x)$. Moreover, we denote
$
(u\cdot\nabla)u
=
\left(
\sum_{j=1}^n u_j\partial_{x_j}u_1,\ldots,
\sum_{j=1}^n u_j\partial_{x_j}u_n
\right)^\top\in \mathbb{R}^n,
$
and $\nabla\cdot u=\sum_{j=1}^n\partial_{x_j}u_j$. We also introduce the cone of nonnegative integrable functions on \(\mathbb R^n\),
$$L^1_+(\mathbb R^n):=
\left\{
\rho\in L^1(\mathbb R^n)
:
\rho\geq 0
\ \text{a.e. in }\mathbb R^n
\right\}.$$

Finally, for $k\in\mathbb N$ and $\alpha\in(0,1)$,  $C^k_b(\mathbb R^n)$ denotes the space of
functions whose derivatives up to order $k$ are bounded and uniformly continuous on
$\mathbb R^n$, endowed with the norm $\sum_{|\beta|\le k}\|D_x^\beta v\|_{L^\infty}$, and by
$C_b^{k+\alpha}(\mathbb R^n)$ its subspace of functions whose derivatives of order $k$ are
$\alpha$-H\"older continuous, with the corresponding norm. \\

We shall consider the following standing assumptions:

\noindent
\paragraph{Initial data:}
\begin{itemize}
\item[(H1)] The initial density satisfies $\rho_0\in L_+^1(\mathbb R^n)\cap L^\infty(\mathbb R^n)$. Moreover, the initial velocity satisfies $u_0\in W^{1,\infty}(\mathbb R^n;\mathbb R^n)$, and $f_0\in W^{1,\infty}(\mathbb R^n)$.
\end{itemize}
\paragraph{Kernels:}
\begin{itemize}
\item[(H2)] The communication kernel is radial defined by the non-negative function $\phi\in W^{1,\infty}([0,+\infty))$.
\item[(H3)] The chemoattractant kernel is defined by $\psi\in W^{1,\infty}(\mathbb R^n)$.
\end{itemize}

\medskip
\noindent
Under \textup{(H2)} we set $\Phi(x):=\phi(|x|)$, $x\in\mathbb R^n$. Since
\[
|\Phi(x)-\Phi(y)|\le\|\phi'\|_{L^\infty}\bigl|\,|x|-|y|\,\bigr|\le\|\phi'\|_{L^\infty}|x-y|,
\qquad x,y\in\mathbb R^n,
\]
the function $\Phi$ belongs to $W^{1,\infty}(\mathbb R^n)$, with
$\|\Phi\|_{L^\infty(\mathbb R^n)}=\|\phi\|_{L^\infty(0,+\infty)}$ and
$\|\nabla\Phi\|_{L^\infty}\le\|\phi'\|_{L^\infty}$, its gradient being defined for almost every
$x\in\mathbb R^n$.

\medskip
\noindent\textbf{Parabolic coefficients.} Set $a(x,t):=\sigma(x,t)\sigma(x,t)^\top$, $(x,t)\in\mathbb R^n\times[0,T]$. We assume that, for some exponent \(\alpha\in(0,1)\), the following conditions hold:

\begin{itemize}

\item[\textup{(C1)}]
The matrix-valued function $a$ and the vector field $b$ are bounded and parabolically Hölder continuous. More precisely, there exists $L>0$ such that
\[
\begin{aligned}
&\|a(x,t)\|+|b(x,t)|\le L,\
&\|a(x,t)-a(y,s)\|
+
|b(x,t)-b(y,s)|
\le
L\left(|x-y|^\alpha+|t-s|^{\alpha/2}\right)
\end{aligned}
\]
for every \(x,\,y\in\mathbb R^n\) and \(t,\,s\in[0,T]\).

\item[\textup{(C2)}]  There exists \(\mu_0>0\) such that $
\left\langle a(x,t)\xi,\xi\right\rangle
\ge \mu_0|\xi|^2$, for every \(x,\xi\in\mathbb R^n\) and \(t\in[0,T]\).

\end{itemize}

\noindent 
We recall that under assumptions \textup{(C1)-(C2)}, the operator $\mathcal P:=\partial_t-
{\rm tr}\!\left(a(x,t)\nabla^2\right)+b(x,t)\cdot\nabla$ admits a fundamental solution $G=G(x,t;\xi,s)$, with $x,\,\xi\in\mathbb R^n$, and $0\le s<t\le T$ . Moreover, for every multi-index \(\beta\) with \(0\le|\beta|\le2\),
there exist constants \(C_T,\,c_T>0\), depending only $
n,\ T,\ \alpha,\ \mu_0, \ L$, such that
\begin{equation}
\left | D_x^\beta G(x,t;\xi,s)\right|
\le
\frac{C_T}{(t-s)^{(n+|\beta|)/2}}
\exp\left(
-c_T\frac{|x-\xi|^2}{t-s}
\right),
\qquad 0\le s<t\le T.
\label{eq:fundamental_solution_gaussian}
\end{equation}
We also recall that $G$ satisfies
\begin{equation}
\int_{\mathbb R^n}G(x,t;\xi,s)\,d\xi=1,
\qquad 0\le s<t\le T,
\label{eq:fundamental_solution_mass}
\end{equation}
and therefore
\begin{equation}
\int_{\mathbb R^n}
D_x^\beta G(x,t;\xi,s) \,d\xi=0,
\qquad 1\le|\beta|\le2.
\label{eq:fundamental_solution_cancellation}
\end{equation}
For further details on the fundamental solutions and relating properties see \cite{friedman2008partial}.

In the considered setting, the main result concerning local well-posedness of system \eqref{eq:model_intro} is given in Theorem \ref{thm:local-wellposedness}. First, we give the following definition of solution to system \eqref{eq:model_intro}.

\begin{Definition}
\label{def:strong-solution}
Let \(T>0\). A triple $(\rho,u,f)$ is called a solution on \([0,T]\) of the system (\ref{eq:model_intro}) if the following properties are satisfied:
\begin{enumerate}
\item $\rho\in
C\bigl([0,T];L^1_+(\mathbb R^n)\bigr)
\cap
L^\infty\bigl(0,T;L^\infty(\mathbb R^n)\bigr)$.

\item $u\in
C\bigl([0,T];L^\infty(\mathbb R^n;\mathbb R^n)\bigr)
\cap
L^\infty\bigl(0,T;W^{1,\infty}(\mathbb R^n;\mathbb R^n)\bigr)$.

\item $f\in
C\bigl([0,T];L^\infty(\mathbb R^n)\bigr)
\cap
L^\infty\bigl(0,T;W^{1,\infty}(\mathbb R^n)\bigr)$, and
$\nabla^2f\in
L^1\bigl(0,T;L^\infty(\mathbb R^n)\bigr)$.

\item The continuity equation in (\ref{eq:model_intro}) holds in the sense of distributions on
\((0,T)\times\mathbb R^n\), namely,
\begin{equation}
\label{eq:weak-continuity-equation}
\int_0^T\int_{\mathbb R^n}
\rho
\left(
\partial_t\varphi
+
u\cdot\nabla\varphi
\right)
\,dx\,dt
+
\int_{\mathbb R^n}
\rho_0(x)\varphi(0,x)\,dx
=
0
\end{equation}
for every
$\varphi\in C_c^1([0,T)\times\mathbb R^n)$.

\item The equation for the velocity field and the parabolic equation in (\ref{eq:model_intro}) hold almost everywhere in
\((0,T)\times\mathbb R^n\).

\item
The initial conditions are attained in the sense that
$\rho(t)\longrightarrow\rho_0$ in $L^1(\mathbb R^n)$,
$u(t)\longrightarrow u_0$ in $L^\infty(\mathbb R^n;\mathbb R^n)$,
$f(t)\longrightarrow f_0$ in $L^\infty(\mathbb R^n)$, and
\[
\nabla u(t)\ \overset{*}{\rightharpoonup}\ \nabla u_0,
\qquad
\nabla f(t)\ \overset{*}{\rightharpoonup}\ \nabla f_0
\qquad
\text{weakly-}*\ \text{in }L^\infty(\mathbb R^n),
\]
as \(t\to0^+\).
\end{enumerate}
\end{Definition}

 We first present two auxilary Lemma, concerning parabolic and alignment estimates respectively.

\begin{lemma}[Parabolic estimates]
\label{lem:parabolic_estimates_matrix}
Let \textup{(C1)-(C2)} hold and let $\rho\in
C\left([0,T];L^1(\mathbb R^n)\right)$ satisfy
\begin{equation}
\sup_{t\in[0,T]}
\|\rho(\cdot,t)\|_{L^1}
\le M_\rho.
\label{eq:rho_mass_bound_parabolic}
\end{equation}
Assume also that $\psi\in W^{1,\infty}(\mathbb R^n)$, $f_0\in W^{1,\infty}(\mathbb R^n)$, and let \(f\) solve
\begin{equation}
\begin{aligned}
\partial_t f - {\rm tr}\left(a(x,t)\nabla^2f\right)
+b(x,t)\cdot\nabla f
&=
(\psi*\rho)(x,t),
&&x\in\mathbb R^n,\quad t\in(0,T]
\end{aligned}
\label{eq:parabolic_problem_lemma}
\end{equation}
with $f(x,0)=f_0(x)$, for every $x\in\mathbb R^n$. Then there exists a constant \(C_T>0\), depending only on \(T\) and
on the constants appearing in \textup{(C1)-(C2)}, such that
\begin{equation}
\|\nabla f(\cdot,t)\|_{L^\infty}
\le
C_T\|\nabla f_0\|_{L^\infty}
+
C_T\|\psi\|_{L^\infty}M_\rho\,t^{1/2},
\qquad 0<t\le T,
\label{eq:parabolic_gradient_estimate}
\end{equation}
and
\begin{equation}
\|\nabla ^2 f(\cdot,t)\|_{L^\infty}
\le
C_T\|\nabla f_0\|_{L^\infty}t^{-1/2}
+
C_T\|\nabla\psi\|_{L^\infty}M_\rho\,t^{1/2},
\qquad 0<t\le T.
\label{eq:parabolic_hessian_estimate}
\end{equation}
In particular, $\nabla^2 f \in L^1\left([0,T];L^\infty(\mathbb R^n)\right)$.
\end{lemma}

\begin{proof}\rm
Let $h(x,t):=(\psi*\rho)(x,t)$. From \eqref{eq:rho_mass_bound_parabolic} follows
$\|h(\cdot,t)\|_{L^\infty}\le\|\psi\|_{L^\infty}M_\rho$, and differentiating with respect to $x$,
\begin{equation}
\|\nabla h(\cdot,t)\|_{L^\infty}
\le
\|\nabla\psi\|_{L^\infty}
\|\rho(\cdot,t)\|_{L^1}
\le
\|\nabla\psi\|_{L^\infty}M_\rho.
\label{eq:source_gradient_bound}
\end{equation}
Using the representation formula associated with the fundamental solution, we write
$f=f^{(0)}+f^{(h)}$, where
\[
f^{(0)}(x,t)
:=
\int_{\mathbb R^n}
G(x,t;\xi,0)f_0(\xi)\,d\xi,\quad
f^{(h)}(x,t)
:=
\int_0^t\int_{\mathbb R^n}
G(x,t;\xi,s)h(\xi,s)\,d\xi ds.
\]
From \eqref{eq:fundamental_solution_cancellation}, we have
$
\nabla f^{(0)}(x,t)=
\int_{\mathbb R^n}
\nabla G(x,t;\xi,0)
\bigl(f_0(\xi)-f_0(x)\bigr)\,d\xi.
$
Since \(f_0\in W^{1,\infty}(\mathbb R^n)\), it holds $|f_0(\xi)-f_0(x)|
\le
|\nabla f_0|_{L^\infty}|x-\xi|$. From \eqref{eq:fundamental_solution_gaussian}, using the change of variable $z=\frac{x-\xi}{\sqrt t}$, we get
\[
\begin{aligned}
\|\nabla f^{(0)}(x,t)\|
&\le
C_T\|\nabla f_0\|_{L^\infty}
t^{-(n+1)/2}
\int_{\mathbb R^n}
|x-\xi|
\exp\left(
-c_T\frac{|x-\xi|^2}{t}
\right)d\xi\\
&=C_T\|\nabla f_0\|_{L^\infty}\int_{\mathbb R^n}|z|e^{-c_T|z|^2}\,dz.
\end{aligned}
\]
Hence
\begin{equation}
\|\nabla f^{(0)}(\cdot,t)\|_{L^\infty}
\le
C_T\|\nabla f_0\|_{L^\infty}.
\label{eq:initial_gradient_bound}
\end{equation}

\medskip
\noindent
Differentiating the representation formula and using
\eqref{eq:fundamental_solution_gaussian}, we find
\[
\begin{aligned}
|\nabla f^{(h)}(x,t)|
&\le C_T
\sup_{0\le s\le t}\|h(\cdot,s)\|_{L^\infty}
\int_0^t
(t-s)^{-(n+1)/2}
&
\int_{\mathbb R^n}
\exp\left(
-c_T\frac{|x-\xi|^2}{t-s}
\right)d\xi\,ds.
\end{aligned}
\]
Since
$\displaystyle \int_{\mathbb R^n}
\exp\left(
-c_T\frac{|x-\xi|^2}{t-s}
\right)d\xi=
C(t-s)^{n/2}$, with $C=C(n,c_T)$, we get 
\[
|\nabla f^{(h)}(x,t)|
\le
C_T
\sup_{0\le s\le t}\|h(\cdot,s)\|_{L^\infty}
\int_0^t(t-s)^{-1/2}\,ds
=
C_T
\sup_{0\le s\le t}\|h(\cdot,s)\|_{L^\infty}
2t^{1/2}\,ds.
\]
By definition of $h$ we hence get 
\begin{equation}
\|\nabla f^{(h)}(\cdot,t)\|_{L^\infty}
\le
C_T\|\psi\|_{L^\infty}M_\rho t^{1/2}.
\label{eq:source_gradient_estimate}
\end{equation}
Combining \eqref{eq:initial_gradient_bound} and
\eqref{eq:source_gradient_estimate} proves
\eqref{eq:parabolic_gradient_estimate}. 
The estimate in \eqref{eq:parabolic_hessian_estimate} follows by the same parabolic regularity argument, recalling that $\nabla(\psi*\rho)=(\nabla\psi)*\rho$ and replacing the estimate for the source term with the corresponding estimate for its gradient. 

\end{proof}

\begin{remark}\label{rem:more_regular_f0}\rm 
We remark that the singular factor \(t^{-1/2}\) in
\eqref{eq:parabolic_hessian_estimate} is due to the minimal assumption $f_0\in W^{1,\infty}(\mathbb R^n)$ considered. If the initial datum has the stronger regularity $f_0\in C_b^{2+\alpha}(\mathbb R^n)$, then the classical Schauder estimates for uniformly parabolic equations under assumptions \textup{(C1)-(C2)} imply
\[
\sup_{0\le t\le T}
\|\nabla ^2 f(\cdot,t)\|_{L^\infty}
\le
C_T\left(
\|f_0\|_{C_b^{2+\alpha}}
+
\|\nabla\psi\|_{L^\infty}M_\rho T^{1/2}
\right).
\]
Thus, under stronger regularity, the Hessian of \(f\) is uniformly bounded up to the initial time. 
\end{remark}

Let $\mathcal{M}_{+}(\mathbb{R}^{n})$ denote the space of finite nonnegative Radon measures on
\(\mathbb{R}^{n}\). For
\(\mu,\nu\in\mathcal{M}_{+}(\mathbb{R}^{n})\), their
bounded-Lipschitz distance is defined by
\begin{equation}
\label{eq:BL-distance}
d_{\mathrm{BL}}(\mu,\nu)
:=
\sup
\left\{
\left|
\int_{\mathbb{R}^{n}}\zeta\,d(\mu-\nu)
\right|
:
\|\zeta\|_{L^\infty}
+
\operatorname{Lip}(\zeta)
\leq 1
\right\},
\end{equation}
Accordingly, we set
\begin{equation}
\label{eq:BL-norm}
\|\zeta\|_{\mathrm{BL}}:=\|\zeta\|_{L^\infty}+\operatorname{Lip}(\zeta),
\end{equation}
so that $d_{\mathrm{BL}}(\mu,\nu)=\sup\bigl\{\bigl|\int_{\mathbb R^n}\zeta\,d(\mu-\nu)\bigr|:
\|\zeta\|_{\mathrm{BL}}\le1\bigr\}$, and we note that
$\|\zeta(x-\cdot)\|_{\mathrm{BL}}=\|\zeta\|_{\mathrm{BL}}$ for every $x\in\mathbb R^n$. 

Whenever \(\rho_1,\rho_2\in L^1_+(\mathbb{R}^n)\), we use the
abbreviated notation
\[
d_{\mathrm{BL}}(\rho_1,\rho_2)
:=
d_{\mathrm{BL}}
\bigl(\rho_1(x)\,dx,\rho_2(x)\,dx\bigr).
\]
For additional details related to this notion of distance see, for instance \cite{dudley2018real} (Chapter 11).

In the following Lemma we establish some 
useful estimates for the alignment operator, defined 
for $\rho\in L^1_+(\mathbb R^n)$ and $u\in L^\infty(\mathbb R^n;\mathbb R^n)$ by
\begin{equation}
\label{eq:alignment-operator}
\mathcal A[\rho,u](x)
:=
\int_{\mathbb R^n}\Phi(x-y)\bigl(u(y)-u(x)\bigr)\rho(y)\,dy,
\qquad x\in\mathbb R^n,
\end{equation}

\begin{lemma}[Alignment estimates]
\label{prop:alignment-estimates}

Under assumption {\rm (H2)}, the following statements hold:
\begin{enumerate}

\item
Let $\rho\in L^1_+(\mathbb R^n)$ and $u\in
W^{1,\infty}(\mathbb R^n;\mathbb R^n)$, then $\mathcal A[\rho,u]
\in
W^{1,\infty}(\mathbb R^n;\mathbb R^n)$ and
\begin{equation}
\label{eq:alignment-strong-estimate}
\|\mathcal A[\rho,u]\|_{W^{1,\infty}}
\leq
C_\Phi
\|\rho\|_{L^1}
\|u\|_{W^{1,\infty}},
\end{equation}
where one may take $C_\Phi=
3\|\Phi\|_{L^\infty}
+
2\|\nabla\Phi\|_{L^\infty}$.
\item
Let $\rho_1,\rho_2\in L^1_+(\mathbb R^n)$ satisfy $\|\rho_i\|_{L^1}\leq M$, for i=1,2, and let $u_1$, $u_2$
$\in W^{1,\infty}(\mathbb R^n;\mathbb R^n)$ satisfy $\|u_i\|_{W^{1,\infty}}\leq R$, for i=1,2, then
\begin{eqnarray}\label{eq:alignment-stability-estimate}
\!\!\!\!\!\!\!\!\|
\mathcal A[\rho_1,u_1]-\mathcal A[\rho_2,u_2]
\|_{L^\infty}
\leq{}
\sqrt n\,C_{\Phi,R}\,
d_{\mathrm{BL}}(\rho_1,\rho_2)
+
2M\|\Phi\|_{L^\infty}\|u_1-u_2\|_{L^\infty},
\end{eqnarray}
where $C_{\Phi,R}=R\left(3\|\Phi\|_{L^\infty}+2\|\nabla\Phi\|_{L^\infty}\right)$.
\end{enumerate}
\end{lemma}

\begin{proof}\rm
\begin{enumerate}
\item  Since $|u(y)-u(x)|
\leq
2\|u\|_{L^\infty}
$, for any $x$, $y$, it follows that
\begin{equation}
\label{eq:alignment-Linf-proof}
\|\mathcal A[\rho,u]\|_{L^\infty}
\leq
2\|\Phi\|_{L^\infty}
\|\rho\|_{L^1}
\|u\|_{L^\infty}.
\end{equation}
For almost every \(x\in\mathbb R^n\), differentiation with respect
to \(x\) gives

\begin{eqnarray}
\!\!\!\!\nabla \mathcal A[\rho,u](x)
=
\int_{\mathbb R^n}
\bigl(u(y)-u(x)\bigr)\otimes\nabla\Phi(x-y)\,
\rho(y)\,dy
-
\left(
\int_{\mathbb R^n}\Phi(x-y)\rho(y)\,dy
\right)
\nabla u(x).
\end{eqnarray}

Consequently,
\begin{eqnarray}
\!\!\!\!\!\!|\nabla\mathcal A[\rho,u](x)|
\leq
2\|\nabla\Phi\|_{L^\infty}
\|\rho\|_{L^1}
\|u\|_{L^\infty}
+
\|\Phi\|_{L^\infty}
\|\rho\|_{L^1}
\|\nabla u\|_{L^\infty}.
\label{eq:alignment-bo}
\end{eqnarray}
Combining \eqref{eq:alignment-bo} with
\eqref{eq:alignment-Linf-proof}, we obtain
$
\|\mathcal A[\rho,u]\|_{W^{1,\infty}}
\leq
\left(
3\|\Phi\|_{L^\infty}
+
2\|\nabla\Phi\|_{L^\infty}
\right)
\|\rho\|_{L^1}
\|u\|_{W^{1,\infty}},
$
which proves \eqref{eq:alignment-strong-estimate}.

\item We write $\mathcal A[\rho_1,u_1]
-
\mathcal A[\rho_2,u_2]
=
\mathcal A[\rho_1,u_1]
-
\mathcal A[\rho_2,u_1]
+
\mathcal A[\rho_2,u_1]
-
\mathcal A[\rho_2,u_2]$.

For any fixed \(x\in\mathbb R^n\) we rewrite
\[
\mathcal A[\rho_1,u_1](x)
-
\mathcal A[\rho_2,u_1](x)
=
\int_{\mathbb R^n}
K_x(y)
\bigl(\rho_1(y)-\rho_2(y)\bigr)\,dy,
\]

with
$
K_x(y)
:=
\Phi(x-y)
\bigl(u_1(y)-u_1(x)\bigr)
$.
Hence $\|K_x\|_{L^\infty}
\leq
2\|\Phi\|_{L^\infty}
\|u_1\|_{L^\infty}$.
Moreover, for any \(y,z\in\mathbb R^n\),
\[
\begin{aligned}
|K_x(y)-K_x(z)|
\leq{}&
|\Phi(x-y)-\Phi(x-z)|
\,|u_1(y)-u_1(x)|
\\
&+
|\Phi(x-z)|
\,|u_1(y)-u_1(z)|.
\end{aligned}
\]
Hence
\[
\operatorname{Lip}(K_x)
\leq
2\|\nabla\Phi\|_{L^\infty}
\|u_1\|_{L^\infty}
+
\|\Phi\|_{L^\infty}
\|\nabla u_1\|_{L^\infty}.
\]
Using $\|u_1\|_{W^{1,\infty}}\leq R$, we infer that
$
\|K_x\|_{L^\infty}
+
\operatorname{Lip}(K_x)
\leq
C_{\Phi,R}$.

Since $|K_x^i|\le|K_x|$ and $|K_x^i(y)-K_x^i(z)|\le|K_x(y)-K_x(z)|$ for every
$i\in\{1,\dots,n\}$, each scalar component $K_x^i$ satisfies
$\|K_x^i\|_{L^\infty}+\operatorname{Lip}(K_x^i)\le C_{\Phi,R}$, and is therefore an admissible
test function for \eqref{eq:BL-distance} up to normalization. Consequently,
\[
\left|\int_{\mathbb R^n}K_x^i(y)\bigl(\rho_1(y)-\rho_2(y)\bigr)\,dy\right|
\le
C_{\Phi,R}\,d_{\mathrm{BL}}(\rho_1,\rho_2),
\qquad i=1,\dots,n,
\]
and, summing the squares of the $n$ components,
\[
\left|\int_{\mathbb R^n}K_x(y)\bigl(\rho_1(y)-\rho_2(y)\bigr)\,dy\right|
\le
\sqrt n\,C_{\Phi,R}\,d_{\mathrm{BL}}(\rho_1,\rho_2).
\]
Taking the supremum with respect to \(x\), we obtain
\begin{equation}
\label{eq:alignment-density-stability}
\|\mathcal A[\rho_1,u_1]-\mathcal A[\rho_2,u_1]\|_{L^\infty}
\le
\sqrt n\,C_{\Phi,R}\,d_{\mathrm{BL}}(\rho_1,\rho_2).
\end{equation}

Since $\mathcal A$ is linear with respect to its second argument, denoting $v:=u_1-u_2$, we get
$
\mathcal A[\rho_2,u_1]-\mathcal A[\rho_2,u_2]=\mathcal A[\rho_2,v].
$
By \eqref{eq:alignment-Linf-proof},
\begin{eqnarray}
\|
\mathcal A[\rho_2,u_1]
-
\mathcal A[\rho_2,u_2]
\|_{L^\infty}
\leq
2\|\Phi\|_{L^\infty}
\|\rho_2\|_{L^1}
\|u_1-u_2\|_{L^\infty}
\leq
2M\|\Phi\|_{L^\infty}
\|u_1-u_2\|_{L^\infty}.
\label{eq:diff_A_u}
\end{eqnarray}
The result in \eqref{eq:alignment-stability-estimate} follows from \eqref{eq:alignment-density-stability} and \eqref{eq:diff_A_u}.
\end{enumerate}
\end{proof}

\begin{theorem}[Local well-posedness]
\label{thm:local-wellposedness}
Assume {\rm (H1)-(H3)} and {\rm (C1)-(C2)}, and let $\gamma>0$. Then there exists
$T_*\in(0,T]$, 
depending only on $n$, $\gamma$, on the constants appearing in {\rm (H2)-(H3)}, {\rm (C1)-(C2)}, and on the quantities $\|\rho_0\|_{L^1}$, $\|\rho_0\|_{L^\infty}$, $\|u_0\|_{W^{1,\infty}}$, $\|\nabla f_0\|_{L^\infty}$, 
such that system \eqref{eq:model_intro}, with initial datum
$(\rho_0,u_0,f_0)$, admits a solution $(\rho,u,f)$ on $[0,T_*]$ in the sense of
Definition~\ref{def:strong-solution}. 

In particular,
\[
\rho\in
C\bigl([0,T_*];L_+^1(\mathbb R^n)\bigr)
\cap
L^\infty\bigl(0,T_*;L^\infty(\mathbb R^n)\bigr),
\]
\[
u\in
C\bigl([0,T_*];L^\infty(\mathbb R^n;\mathbb R^n)\bigr)
\cap
L^\infty\bigl(0,T_*;W^{1,\infty}(\mathbb R^n;\mathbb R^n)\bigr),
\]
\[
f\in
C\bigl([0,T_*];L^\infty(\mathbb R^n)\bigr)
\cap
L^\infty\bigl(0,T_*;W^{1,\infty}(\mathbb R^n)\bigr),
\qquad
\nabla^2f\in
L^1\bigl(0,T_*;L^\infty(\mathbb R^n)\bigr).
\]
Moreover, 
$\|\rho(t)\|_{L^1}=\|\rho_0\|_{L^1}$ for every $t\in[0,T_*]$, and there exists $R>0$,
depending only on the initial data, such that
\begin{equation}
\label{eq:uniform-bound-solution}
\sup_{t\in[0,T_*]}
\Bigl(
\|\rho(t)\|_{L^\infty}
+
\|u(t)\|_{W^{1,\infty}}
\Bigr)
\le R .
\end{equation}
Finally, the solution is unique in the class of Definition~\ref{def:strong-solution}.
\end{theorem}

\begin{proof}\rm

We split the proof in five steps, in order to highlight the main points of the fixed-point argument considered in our approach.

\medskip
\noindent
{\em Step 1: Definition and completeness of the fixed-point
space.} 

Let $R>0$ and $T>0$, to be chosen below. We denote by
$\mathcal B_{R,T}$ the set of all pairs $(\rho,u)$ satisfying
\[
\rho\in
C\bigl([0,T];(L^1_+(\mathbb R^n),d_{\mathrm{BL}})\bigr)
\cap
L^\infty\bigl(0,T;L^\infty(\mathbb R^n)\bigr),
\]
\[
u\in
C\bigl([0,T];L^\infty(\mathbb R^n;\mathbb R^n)\bigr)
\cap
L^\infty\bigl(0,T;
W^{1,\infty}(\mathbb R^n;\mathbb R^n)\bigr),
\]
together with
$\rho(0)=\rho_0,
u(0)=u_0,$
and
\begin{equation}
\label{eq:fixed-point-mass}
\|\rho(t)\|_{L^1}=m_0
\qquad
\text{for every }t\in[0,T],
\end{equation}
with $m_0:=\|\rho_0\|_{L^1}$.

We also impose the uniform bound
\begin{equation}
\label{eq:fixed-point-ball-bound}
\sup_{t\in[0,T]}
\left(
\|\rho(t)\|_{L^\infty}
+
\|u(t)\|_{W^{1,\infty}}
\right)
\leq R.
\end{equation}
Here, continuity of $\rho$ is understood with respect to the
bounded-Lipschitz distance in \eqref{eq:BL-distance}.
We endow $\mathcal B_{R,T}$ with the metric
\begin{equation}
\label{eq:fixed-point-metric}
{\cal D}_T\bigl((\rho_1,u_1),(\rho_2,u_2)\bigr)
:=
\sup_{t\in[0,T]}
\left[
d_{\mathrm{BL}}\bigl(\rho_1(t),\rho_2(t)\bigr)
+
\|u_1(t)-u_2(t)\|_{L^\infty}
\right].
\end{equation}
We first verify that $(\mathcal B_{R,T},{\cal D}_T)$ is complete. Let $\{(\rho_k,u_k)\}_k$ be a ${\cal D}_T$-Cauchy sequence. Then
$\rho_k$ converges uniformly in time, with respect to
$d_{\mathrm{BL}}$, to a continuous map $t\longmapsto\mu_t$
of finite nonnegative Radon measures.
Moreover, taking $\zeta\equiv1$ in \eqref{eq:BL-distance}, and using
\eqref{eq:fixed-point-mass}, we obtain $\mu_t(\mathbb R^n)=m_0$ for every $t\in[0,T]$.

It remains to show that $\mu_t$ is absolutely continuous with
$L^\infty$ density. For every nonnegative
$\zeta\in C_c(\mathbb R^n)$, from
\eqref{eq:fixed-point-ball-bound} we get
$
\int_{\mathbb R^n}\zeta(x)\rho_k(t,x)\,dx
\leq
R\int_{\mathbb R^n}\zeta(x)\,dx.
$
Passing to the limit in the bounded-Lipschitz topology yields
\[
\int_{\mathbb R^n}\zeta\,d\mu_t
\leq
R\int_{\mathbb R^n}\zeta(x)\,dx.
\]
Hence $\mu_t$ is absolutely continuous with respect to the Lebesgue measure and its Radon-Nikodym density, denoted by $\rho(t,\cdot)$, satisfies $0\leq\rho(t,x)\leq R$, for a.e. $x\in\mathbb R^n$. Thus the density component of the limit remains in
$L^1_+\cap L^\infty$ and retains mass $m_0$.\\
 On the other hand, $u_k$ converges uniformly in
$[0,T]\times\mathbb R^n$ to some limit function $u$. Since $\|\nabla u_k(t)\|_{L^\infty}\leq R$, the functions $u_k(t,\cdot)$ are uniformly Lipschitz. Passing to the
uniform limit gives
$
|u(t,x)-u(t,y)|
\leq
R|x-y|.
$
Hence $u(t)\in W^{1,\infty}(\mathbb R^n;\mathbb R^n)$, $\|\nabla u(t)\|_{L^\infty}\leq R$.\\
\noindent 
Repeating the arguments above along a subsequence realizing the corresponding lower limits gives
$\|\rho(t)\|_{L^\infty}+\|u(t)\|_{W^{1,\infty}}\le
\liminf_k\bigl(\|\rho_k(t)\|_{L^\infty}+\|u_k(t)\|_{W^{1,\infty}}\bigr)\le R$, proving that the constraint \eqref{eq:fixed-point-ball-bound} passes to the limit.
Therefore, the limit belongs to $\mathcal B_{R,T}$, proving completeness of the space.\\

\noindent
From now on we choose $R$ sufficiently large so that
$\mathcal B_{R,T}$ is nonempty and, more precisely,
\begin{equation}
\label{eq:R-choice}
R>
2\left(
\|\rho_0\|_{L^\infty}
+
\|u_0\|_{W^{1,\infty}}
\right)+1.
\end{equation}

\medskip
\noindent
{\em Step 2: Definition of the fixed-point map and invariance of the ball.} Let $(\bar\rho,\bar u)\in\mathcal B_{R,T}$, we define $(\rho,u)=\mathcal T(\bar\rho,\bar u)$, as follows. First, let $f$ be the solution of the linear parabolic problem
\begin{equation}
\label{eq:fixed-f-revised}
\left\{
\begin{aligned}
&\partial_tf
-\operatorname{tr}
 \bigl(\sigma\sigma^\top\nabla^2f\bigr)
+b\cdot\nabla f
=\psi*\bar\rho,
\\
&f(0,\cdot)=f_0.
\end{aligned}
\right.
\end{equation}
Second, let $\rho$ solve the continuity equation
\begin{equation}
\label{eq:fixed-rho-revised}
\left\{
\begin{aligned}
&\partial_t\rho+\nabla\cdot(\rho\bar u)=0,
\\
&\rho(0,\cdot)=\rho_0.
\end{aligned}
\right.
\end{equation}
Finally, let $u$ solve the linear transport equation
\begin{equation}
\label{eq:fixed-u-revised}
\left\{
\begin{aligned}
&\partial_tu+\bar u\cdot\nabla u
=
\mathcal A[\bar\rho,\bar u]+\gamma\nabla f,
\\
&u(0,\cdot)=u_0.
\end{aligned}
\right.
\end{equation}

\noindent
Since $\bar u$ is bounded and uniformly Lipschitz in space, both transport equations are
well posed by the classical theory of characteristics. As for the parabolic problem, the
representation formula associated with the fundamental solution $G$ defines  $f$ which solves \eqref{eq:fixed-f-revised} almost everywhere in $(0,T)\times\mathbb R^n$ and
satisfies $\|f(t)\|_{L^\infty}\le\|f_0\|_{L^\infty}+m_0\|\psi\|_{L^\infty}t$, together with
$f\in C\bigl([0,T];L^\infty(\mathbb R^n)\bigr)$. Moreover, by
Lemma~\ref{lem:parabolic_estimates_matrix}, assumptions {\rm (C1)-(C2)} and
$\|\psi*\bar\rho\|_{L^\infty(0,T;W^{1,\infty})}\le m_0\|\psi\|_{W^{1,\infty}}$, we obtain
\[
f\in C\bigl([0,T];L^\infty(\mathbb R^n)\bigr)\cap L^\infty\bigl(0,T;W^{1,\infty}(\mathbb R^n)\bigr),
\qquad
\nabla^2f\in L^1\bigl(0,T;L^\infty(\mathbb R^n)\bigr).
\]

\medskip
\noindent
First we show that $\mathcal T$ maps $\mathcal B_{R,T}$ into itself.

Let $X$ be the flow generated by $\bar u$, i.e. 
\begin{equation}
\left\{
\begin{aligned}
&\frac{d}{dt}X(t,a)
=
\bar u\bigl(t,X(t,a)\bigr),
\\
&X(0,a)=a.
\label{eq:caratteristiche}
\end{aligned}
\right.
\end{equation}
Since $\|\bar u(t)\|_{W^{1,\infty}}\leq R$, the flow is globally defined on $[0,T]$ and is bi-Lipschitz in the
spatial variable. In particular,
$
|X(t,a)-X(t,\tilde a)|
\leq
e^{Rt}|a-\tilde a|.
$
The solution of the continuity equation is the push-forward of the initial measure \(\rho_0(x)\,dx\) through the flow \(X(t,\cdot)\), namely $\rho(t,x)\,dx = X(t,\cdot)_{\#}\bigl(\rho_0(x)\,dx\bigr), $ where \(X(t,\cdot)_{\#}\mu\) denotes the push-forward of a measure \(\mu\), defined by $$(X(t,\cdot)_{\#}\mu)(B) =\mu\!\left(X(t,\cdot)^{-1}(B)\right),\qquad\mbox{\rm  for every Borel set $B\subset\mathbb R^n$.}$$

\noindent
Consequently, $\rho(t,x)\geq0$, for a.e. $x$, and $\|\rho(t)\|_{L^1}
=
\|\rho_0\|_{L^1}
=
m_0$. The characteristic representation reads
\begin{equation}
\label{eq:rho-characteristic-proof}
\rho\bigl(t,X(t,a)\bigr)
=
\rho_0(a)
\exp\left(
-\int_0^t
\nabla\cdot\bar u
\bigl(s,X(s,a)\bigr)\,ds
\right).
\end{equation}

Since $\|\nabla\cdot \bar u\|_{L^\infty}\leq n\|\nabla\bar u\|_{L^\infty}$ and
$\|\nabla\bar u\|_{L^\infty}\leq R$, we obtain
\begin{equation}
\label{eq:rho-Linfty-revised}
\|\rho(t)\|_{L^\infty}\leq\|\rho_0\|_{L^\infty}e^{nRt}.
\end{equation}
We observe that this estimate uses no spatial derivative of $\rho_0$. Then the map $t\mapsto\rho(t)$ is continuous with respect to
$d_{\mathrm{BL}}$. Indeed, for every bounded Lipschitz function
$\zeta$, it holds
\[
\int_{\mathbb R^n}\zeta(x)\rho(t,x)\,dx
=
\int_{\mathbb R^n}
\zeta\bigl(X(t,a)\bigr)\rho_0(a)\,da,
\]
and the continuity of $X(t,a)$ with respect to $t$, together with
dominated convergence, proves the claim.

\medskip
\noindent
As for the velocity, setting $F=\mathcal A[\bar\rho,\bar u]+\gamma\nabla f$, then by Lemma ~\ref{prop:alignment-estimates} and
Lemma~\ref{lem:parabolic_estimates_matrix}, we deduce that 
\begin{equation}
\label{eq:stimeF}   
F \in L^\infty \bigl(0,T;L^\infty(\mathbb R^n;\mathbb R^n)\bigr)  \qquad \nabla F\in L^1\bigl(0,T;L^\infty(\mathbb R^n;\mathbb R^{n\times n})\bigr), 
\end{equation}
with norms bounded only in terms of the initial data,
\(R\), \(m_0\), \(\gamma\), and the coefficients of the system.\\

\noindent
Let \(X\) be the flow generated by \(\bar u\).
Along the characteristics of \(X\),
equation \eqref{eq:fixed-u-revised} becomes
\[
\frac{d}{dt}
u(t,X(t,a))
=
F(t,X(t,a)).
\]
Hence
\begin{equation}
\label{eq:stimau}
\|u(t)\|_{L^\infty}
\le
\|u_0\|_{L^\infty}
+
\int_0^t
\|F(s)\|_{L^\infty}\,ds.
\end{equation}
In order to estimate the spatial gradient of \(u\), we shall argue on its component. For \(h\neq0\) and \(j\in\{1,\ldots,n\}\), set
\[
\delta_h^j u(t,x)
:=
\frac{u(t,x+he_j)-u(t,x)}{h}.
\]
Applying \(\delta_h^j\) to the transport equation gives
$
\partial_t\delta_h^ju
+
\bar u\cdot\nabla\delta_h^ju
=
\delta_h^jF
-
\delta_h^j\bar u\cdot
\nabla u(\,\cdot+he_j).
$
Since $\|\delta_h^j\bar u(t)\|_{L^\infty}
\leq
\|\nabla\bar u(t)\|_{L^\infty}
\leq R$ and $\|\delta_h^jF(t)\|_{L^\infty}
\leq
\|\partial_jF(t)\|_{L^\infty}$, the \(L^\infty\)-estimate for linear transport equations yields
\[
\|\delta_h^ju(t)\|_{L^\infty}
\leq
\|\delta_h^ju_0\|_{L^\infty}
+
\int_0^t
\left(
R\|\nabla u(s)\|_{L^\infty}
+
\|\partial_jF(s)\|_{L^\infty}
\right)ds.
\]
Passing to the supremum with respect to \(h\) and \(j\), we obtain
$$
\|\nabla u(t)\|_{L^\infty}
\leq
\|\nabla u_0\|_{L^\infty}
+
\int_0^t
\left(
R\|\nabla u(s)\|_{L^\infty}
+
\|\nabla F(s)\|_{L^\infty}
\right)ds.
$$
and, by Grönwall's inequality,
\begin{equation}
\label{eq:stimagradu}
\|\nabla u(t)\|_{L^\infty}
\leq
e^{Rt}
\left[
\|\nabla u_0\|_{L^\infty}
+
\int_0^t
\|\nabla F(s)\|_{L^\infty}\,ds
\right].
\end{equation}
Combining \eqref{eq:stimau} and \eqref{eq:stimagradu}, we obtain
\[
\|u(t)\|_{W^{1,\infty}}
\le
e^{Rt}
\left(
\|u_0\|_{W^{1,\infty}}
+
\omega_u(T)
\right),
\qquad
\omega_u(T):=\int_{0}^{T}\|F(s)\|_{W^{1,\infty}}\,ds ,
\]

where $\omega_u(T)\longrightarrow0$ as $T\downarrow0$, uniformly with respect to
$(\bar\rho,\bar u)\in\mathcal B_{R,T}$ by \eqref{eq:stimeF}.   
\noindent
To prove the constraint \eqref{eq:fixed-point-ball-bound}, from \eqref{eq:rho-Linfty-revised} and the estimate above we get
\begin{equation}
\label{eq:Lambda-def}
\sup_{t\in[0,T]}
\Bigl(\|\rho(t)\|_{L^\infty}+\|u(t)\|_{W^{1,\infty}}\Bigr)
\le
\Lambda(T)
:=
e^{nRT}\|\rho_0\|_{L^\infty}
+
e^{RT}\bigl(\|u_0\|_{W^{1,\infty}}+\omega_u(T)\bigr).
\end{equation}
The function $\Lambda$ is continuous and
$\Lambda(T)\longrightarrow\|\rho_0\|_{L^\infty}+\|u_0\|_{W^{1,\infty}}$ as $T\downarrow0$, a
quantity which is smaller than $R/2$ by the choice \eqref{eq:R-choice} of $R$. 
Hence $\Lambda(T)\le R$ for every $T$ small enough, say $T\le T_1$. 
Since moreover
$\rho(0)=\rho_0$, $u(0)=u_0$ and $\|\rho(t)\|_{L^1}=m_0$, we conclude that
\begin{equation}
\label{eq:invariance-ball}
\mathcal T\bigl(\mathcal B_{R,T}\bigr)\subseteq\mathcal B_{R,T},
\qquad T\in(0,T_1] .
\end{equation}

\medskip
\noindent

{\em Step 3: Contraction property}

Let $(\rho_i,u_i)
=
\mathcal T(\bar\rho_i,\bar u_i)$, $i=1,2$, and set
\[
\delta_\rho
:=
\sup_{t\in[0,T]}
d_{\mathrm{BL}}
\bigl(\bar\rho_1(t),\bar\rho_2(t)\bigr),\qquad \delta_u
:=
\sup_{t\in[0,T]}
\|\bar u_1(t)-\bar u_2(t)\|_{L^\infty}.
\]
Let $f_i$ be the solutions of \eqref{eq:fixed-f-revised}
corresponding to $\bar\rho_i$, and set $g:=f_1-f_2$. Then
\begin{equation}
\label{eq:g-parabolic}
\left\{
\begin{aligned}
&\partial_tg
-\operatorname{tr}
 \bigl(\sigma\sigma^\top\nabla^2g\bigr)
+b\cdot\nabla g
=
\psi*(\bar\rho_1-\bar\rho_2),
\\
&g(0,\cdot)=0.
\end{aligned}
\right.
\end{equation}

For every $x\in\mathbb R^n$ 
$\|\psi(x-\cdot)\|_{\mathrm{BL}}=\|\psi\|_{\mathrm{BL}}=:C_\psi$. 
Therefore,
\begin{equation}
\label{eq:source-BL-estimate}
\| \psi*(\bar\rho_1-\bar\rho_2)(t)
\|_{L^\infty}
\leq
C_\psi
d_{\mathrm{BL}}
\bigl(\bar\rho_1(t),\bar\rho_2(t)\bigr)
\leq
C_\psi\delta_\rho.
\end{equation}
Since $g(0)=0$, estimates \eqref{eq:parabolic_gradient_estimate} and
\eqref{eq:source-BL-estimate} give
\begin{equation}
\label{eq:f-difference-revised}
\|\nabla g(t)\|_{L^\infty}
\leq
C_T t^{1/2}\delta_\rho,
\qquad
t\in[0,T].
\end{equation}

\medskip
\noindent
Let $X_i$ be the flows generated by $\bar u_i$. Since
$
\rho_i(t,x)\,dx
=
X_i(t,\cdot)_{\#}
\bigl(\rho_0(x)\,dx\bigr),
$
for every $\zeta$ satisfying $\|\zeta\|_{\mathrm{BL}}\leq1$, we have
\[
\begin{aligned}
&
\left|
\int_{\mathbb R^n}
\zeta(x)
\bigl(\rho_1(t,x)-\rho_2(t,x)\bigr)\,dx
\right|=
\left|
\int_{\mathbb R^n}
\left[
\zeta\bigl(X_1(t,a)\bigr)
-
\zeta\bigl(X_2(t,a)\bigr)
\right]\rho_0(a)\,da
\right|
\\
&\qquad\leq
\int_{\mathbb R^n}
|X_1(t,a)-X_2(t,a)|
\rho_0(a)\,da \leq
m_0
\sup_{a\in\mathbb R^n}
|X_1(t,a)-X_2(t,a)|.
\end{aligned}
\]

Fix $a\in\mathbb R^n$ and set $Y(t):=X_1(t,a)-X_2(t,a)$. Since $X_1(\cdot,a)$ and $X_2(\cdot,a)$
are absolutely continuous, so is $Y(\cdot)$, with
$\dot Y(t)=\bar u_1\bigl(t,X_1(t,a)\bigr)-\bar u_2\bigl(t,X_2(t,a)\bigr)$ for almost every $t$.
For $\varepsilon>0$ let $\theta_\varepsilon(z):=\sqrt{|z|^2+\varepsilon^2}$, $z\in\mathbb R^n$.
Since $\theta_\varepsilon\in C^1(\mathbb R^n)$ with
$|\nabla\theta_\varepsilon|\le1$, the map $t\mapsto\theta_\varepsilon\bigl(Y(t)\bigr)$ is
absolutely continuous and, for almost every $t$,
\[
\begin{aligned}
\frac{d}{dt}\theta_\varepsilon\bigl(Y(t)\bigr)
&\leq
\left|
\bar u_1\bigl(t,X_1(t,a)\bigr)-\bar u_1\bigl(t,X_2(t,a)\bigr)
\right|
+
\left|
\bar u_1\bigl(t,X_2(t,a)\bigr)-\bar u_2\bigl(t,X_2(t,a)\bigr)
\right|
\\
&\leq
R\,|Y(t)|+\delta_u
\leq
R\,\theta_\varepsilon\bigl(Y(t)\bigr)+\delta_u .
\end{aligned}
\]
Since $\theta_\varepsilon\bigl(Y(0)\bigr)=\varepsilon$, Grönwall's inequality gives
$\theta_\varepsilon\bigl(Y(t)\bigr)\le\varepsilon e^{Rt}+\delta_u\int_0^te^{R(t-s)}\,ds$, and
letting $\varepsilon\to0^+$ we obtain
\begin{equation}
\label{eq:flow-comparison}
|X_1(t,a)-X_2(t,a)|
\leq
\delta_u\int_0^te^{R(t-s)}\,ds
\leq
te^{Rt}\delta_u,
\qquad
t\in[0,T],\ a\in\mathbb R^n .
\end{equation}

Taking the supremum over $a$ and then over the admissible test
functions $\zeta$, we obtain
\begin{equation}
\label{eq:rho-contraction-revised}
d_{\mathrm{BL}}
\bigl(\rho_1(t),\rho_2(t)\bigr)
\leq
m_0te^{Rt}\delta_u.
\end{equation}
Consequently,
\begin{equation}
\label{eq:rho-contraction-sup}
\sup_{t\in[0,T]}
d_{\mathrm{BL}}
\bigl(\rho_1(t),\rho_2(t)\bigr)
\leq
m_0Te^{RT}\delta_u.
\end{equation}

\medskip
\noindent
Finally, set $v:=u_1-u_2$. Subtracting the equations satisfied by $u_1$ and $u_2$ gives
\begin{equation}
\label{eq:v-equation}
\left\{
\begin{aligned}
\partial_tv+\bar u_1\cdot\nabla v
={}&
-(\bar u_1-\bar u_2)\cdot\nabla u_2
+
\mathcal A[\bar\rho_1,\bar u_1]
-
\mathcal A[\bar\rho_2,\bar u_2]
+
\gamma\nabla g,
\\
v(0,\cdot)=0\;\;\;\qquad&\qquad\qquad\qquad\qquad.
\end{aligned}
\right.
\end{equation}
Since $\mathcal T$ maps $\mathcal B_{R,T}$ into itself,
$\|\nabla u_2(t)\|_{L^\infty}\leq R$. It follows that
\begin{equation}
\label{eq:transport-difference-term}
\|
(\bar u_1-\bar u_2)\cdot\nabla u_2
\|_{L^\infty}
\leq
R\delta_u.
\end{equation}
By the estimate in
Lemma~\ref{prop:alignment-estimates},
\begin{equation}
\label{eq:alignment-difference-revised}
\left\|
\mathcal A[\bar\rho_1,\bar u_1]
-
\mathcal A[\bar\rho_2,\bar u_2]
\right\|_{L^\infty}
\leq
C_{R,m_0,\phi}
\left(
\delta_\rho+\delta_u
\right).
\end{equation}

\noindent
Integrating \eqref{eq:v-equation} along the characteristics of
$\bar u_1$ and using
\eqref{eq:f-difference-revised},
\eqref{eq:transport-difference-term} and
\eqref{eq:alignment-difference-revised}, we obtain
\[
\begin{aligned}
\|v(t)\|_{L^\infty}
&\leq
C_R
\int_0^t
(\delta_\rho+\delta_u)\,ds
+
C\gamma
\int_0^t
s^{1/2}\delta_\rho\,ds
&\leq
C_{R,\gamma}
\bigl(t+t^{3/2}\bigr)
(\delta_\rho+\delta_u).
\end{aligned}
\]
Therefore,
\begin{equation}
\label{eq:u-contraction-revised}
\sup_{t\in[0,T]}
\|u_1(t)-u_2(t)\|_{L^\infty}
\leq
C_{R,\gamma}
\bigl(T+T^{3/2}\bigr)
(\delta_\rho+\delta_u).
\end{equation}
Combining \eqref{eq:rho-contraction-sup} and
\eqref{eq:u-contraction-revised}, we conclude that
\[
\begin{aligned}
&{\cal D}_T\bigl(
\mathcal T(\bar\rho_1,\bar u_1),
\mathcal T(\bar\rho_2,\bar u_2)
\bigr)\leq
q(T)
{\cal D}_T\bigl(
(\bar\rho_1,\bar u_1),
(\bar\rho_2,\bar u_2)
\bigr),
\end{aligned}
\]
where
\begin{equation}
\label{eq:contraction-factor}
q(T)
:=
m_0Te^{RT}
+
2C_{R,\gamma}\bigl(T+T^{3/2}\bigr).
\end{equation}

Since $q(T)\longrightarrow0$ as $T\downarrow0$, we may choose $T_*\in(0,T_1]$ such that
$q(T_*)<1$. 
By \eqref{eq:invariance-ball}, $\mathcal T$ maps 
$(\mathcal B_{R,T_*},\mathcal D_{T_*})$ into itself and is a contraction there. Therefore,
Banach's fixed-point theorem provides a unique fixed point $(\rho,u)\in\mathcal B_{R,T_*}$.

We here observe, for later use, that both the time $T_1$ in \eqref{eq:invariance-ball} and the time
$T_*$ chosen above depend only on $R$, determined by \eqref{eq:R-choice}, on $m_0$, $n$, $\gamma$,
$\|\nabla f_0\|_{L^\infty}$ and on the constants appearing in {\rm (H2)-(H3)} and {\rm (C1)-(C2)}.
Since the right-hand sides of \eqref{eq:rho-Linfty-revised}, $\Lambda$ in \eqref{eq:Lambda-def}, and $q$ in \eqref{eq:contraction-factor}, are nondecreasing with respect to each of these quantities, the time $T_*$ may be chosen as a nonincreasing function of them. In particular, a uniform bound on the data yields a uniform existence time.

\medskip
\noindent
{\em Step 4: Strong $L^1$ continuity of the density.} 
The fixed-point construction initially
provides continuity of $t\mapsto\rho(t)$ with respect to $d_{\mathrm{BL}}$. We now prove the
stronger property $\rho\in C\bigl([0,T_*];L^1(\mathbb R^n)\bigr)$ required in Definition \ref{def:strong-solution}\\

\noindent
For $0\le s\le t\le T_*$, let $X_{s,t}$ denote the two-parameter flow generated by $u$, that is
$\frac{d}{dt}X_{s,t}(a)=u\bigl(t,X_{s,t}(a)\bigr)$, $X_{s,s}(a)=a$, and let
$A_{s,t}:=X_{s,t}^{-1}$. By the semigroup property of the flow,
$\rho(t)=\bigl(X_{s,t}\bigr)_{\#}\rho(s)$. Since $\|u\|_{L^\infty}\le R$ and
$\|\nabla\cdot u\|_{L^\infty}\le nR$, we have the two uniform bounds
\begin{equation}
\label{eq:two-parameter-flow-bounds}
\|A_{s,t}-{\rm id}\|_{L^\infty}\le R|t-s|,
\qquad
\bigl\|\det\nabla A_{s,t}-1\bigr\|_{L^\infty}\le e^{nR|t-s|}-1 .
\end{equation}
Let $g\in C_c(\mathbb R^n)$. The push-forward $\bigl(X_{s,t}\bigr)_{\#}(g\,dx)$ has density
$g\bigl(A_{s,t}(x)\bigr)\det\nabla A_{s,t}(x)$, so that
\[
\bigl\|(X_{s,t})_{\#}g-g\bigr\|_{L^1}
\le
\bigl\|\det\nabla A_{s,t}-1\bigr\|_{L^\infty}
\bigl\|g\circ A_{s,t}\bigr\|_{L^1}
+
\bigl\|g\circ A_{s,t}-g\bigr\|_{L^1}.
\]
The first term vanishes as $t\to s$ by \eqref{eq:two-parameter-flow-bounds} and the change of
variables formula, which gives
$\|g\circ A_{s,t}\|_{L^1}\le e^{nR|t-s|}\|g\|_{L^1}$; the second one vanishes by the uniform
continuity of $g$, its compact support, and the first bound in
\eqref{eq:two-parameter-flow-bounds}.\\

\noindent
Let now $s\in[0,T_*]$ and $\varepsilon>0$ be fixed, and choose $g\in C_c(\mathbb R^n)$ with
$\|\rho(s)-g\|_{L^1}\le\varepsilon$. Since the push-forward through a bi-Lipschitz bijection
preserves the total variation of signed measures,
$\|(X_{s,t})_{\#}\bigl(\rho(s)-g\bigr)\|_{L^1}=\|\rho(s)-g\|_{L^1}$, and therefore
\[
\|\rho(t)-\rho(s)\|_{L^1}
\le
\bigl\|(X_{s,t})_{\#}\bigl(\rho(s)-g\bigr)\bigr\|_{L^1}
+
\bigl\|(X_{s,t})_{\#}g-g\bigr\|_{L^1}
+
\|g-\rho(s)\|_{L^1}
\le
2\varepsilon+\bigl\|(X_{s,t})_{\#}g-g\bigr\|_{L^1}.
\]
Letting $t\to s$ and then $\varepsilon\to0^+$ proves that $\rho\in C\bigl([0,T_*];L^1(\mathbb R^n)\bigr)$.\\

Finally, we here observe that the initial data are attained in the sense of Definition~\ref{def:strong-solution} (point 6). Indeed the strong convergences follow from the continuity in time established above, while for the gradients
it suffices to integrate by parts against $\varphi\in C_c^\infty(\mathbb R^n)$ and to use the
density of $C_c^\infty(\mathbb R^n)$ in $L^1(\mathbb R^n)$, together with the bounds
$\|\nabla u(t)\|_{L^\infty}\le R$, 
and
$\|\nabla f(t)\|_{L^\infty}\le C_T\bigl(\|\nabla f_0\|_{L^\infty}+\|\psi\|_{L^\infty}m_0T_*^{1/2}\bigr)$, which follows from \eqref{eq:parabolic_gradient_estimate}.\\

\medskip
\noindent
{\em Step 5: Uniqueness in the solution class.}  Banach’s
theorem ensures the uniqueness of the fixed point in the ball \(\mathcal B_{R,T_*}\). In the following we prove uniqueness in the larger class of solutions of Definition~\ref{def:strong-solution}, concluding the proof.

\noindent Let $(\rho_1,u_1,f_1)$,  $(\rho_2,u_2,f_2)$ be two solutions on \([0,T]\), in the sense of
Definition~\ref{def:strong-solution}, with the same initial data. We denote
\[
\delta u:=u_1-u_2,
\qquad
\delta f:=f_1-f_2.
\]
Since both solutions belong to the class of
Definition~\ref{def:strong-solution}, there exists a constant
\(M>0\) such that, for \(i=1,2\),
\[
\sup_{t\in[0,T]}
\left(
 \|\rho_i(t)\|_{L^\infty}
 +
 \|u_i(t)\|_{W^{1,\infty}}
 +
 \|f_i(t)\|_{W^{1,\infty}}
\right)
\leq M.
\]
Moreover, mass conservation gives
$
\|\rho_i(t)\|_{L^1}
=
\|\rho_0\|_{L^1}
=:m_0$, $t\in[0,T].$

For \(0<\tau\leq T\), we define
\[
\mathcal R_\tau
:=
\sup_{t\in[0,\tau]}
d_{\mathrm{BL}}\bigl(\rho_1(t,x)\,dx,\rho_2(t,x)\,dx\bigr),\qquad 
\mathcal U_\tau
:=
\sup_{t\in[0,\tau]}
\|\delta u(t)\|_{L^\infty}.
\]

We first estimate the difference between the two densities. 
Since $u_i\in L^\infty\bigl(0,T;W^{1,\infty}\bigr)\cap C\bigl([0,T];L^\infty\bigr)$, the vector
field $u_i$ is bounded and Lipschitz continuous in space, uniformly in time. Hence the flows
$X_i$ generated by \(u_i\) are well defined, bi-Lipschitz in the spatial variable, and the Cauchy problem for the
continuity equation with velocity $u_i$ and datum $\rho_0$ admits a unique distributional
solution in $C\bigl([0,T];L^1(\mathbb R^n)\bigr)$, which is given by the push-forward of the
initial measure along the flow, namely
$\rho_i(t,x)\,dx=X_i(t,\cdot)_{\#}\bigl(\rho_0(x)\,dx\bigr)$, see \cite{ambrosio2008gradient} (Chapter~8).

The comparison of the two flows is obtained exactly as in \eqref{eq:flow-comparison}, replacing $\bar u_i$ with $u_i$
and $R$ with $M$ as the common Lipschitz constant. We obtain that for every
$a\in\mathbb R^n$ and every $t\in[0,\tau]$,
\[
|X_1(t,a)-X_2(t,a)|
\leq
\int_0^t e^{M(t-s)}\|\delta u(s)\|_{L^\infty}\,ds
\leq
\tau e^{M\tau}\,\mathcal U_\tau .
\]

Hence, for every test function \(\zeta\) satisfying
\(\|\zeta\|_{\mathrm{BL}}\leq1\), we have
\[
\begin{aligned}
\left|
\int_{\mathbb R^n}
\zeta(x)\bigl(\rho_1(t,x)-\rho_2(t,x)\bigr)\,dx
\right|
&=
\left|
\int_{\mathbb R^n}
\bigl[
\zeta(X_1(t,a))-\zeta(X_2(t,a))
\bigr]
\rho_0(a)\,da
\right|
\\
&\leq
m_0
\sup_{a\in\mathbb R^n}
|X_1(t,a)-X_2(t,a)|
\leq
m_0 t e^{Mt}\mathcal U_\tau .
\end{aligned}
\]
Taking the supremum over the admissible test functions and over
\(t\in[0,\tau]\), we obtain
\begin{equation}
\label{eq:uniqueness-density-estimate}
\mathcal R_\tau
\leq
m_0\tau e^{M\tau}\mathcal U_\tau.
\end{equation}
We next estimate the difference between the parabolic components.
The function $\delta f:=f_1-f_2$ solves
\[
\partial_t\delta f
-
\operatorname{tr}
\bigl(\sigma\sigma^\top\nabla^2\delta f\bigr)
+
b\cdot\nabla\delta f
=
\psi*(\rho_1-\rho_2),
\qquad
\delta f(0)=0.
\]
For every \(x\in\mathbb R^n\) and \(t\in[0,\tau]\), the definition of
the bounded-Lipschitz distance gives
\[
\begin{aligned}
\left|
\psi*(\rho_1-\rho_2)(t,x)
\right|
&\leq
\|\psi(x-\cdot)\|_{\mathrm{BL}}\,
d_{\mathrm{BL}}\bigl(\rho_1(t)\,dx,\rho_2(t)\,dx\bigr)
&\leq
\|\psi\|_{\mathrm{BL}}\mathcal R_\tau.
\end{aligned}
\]

\noindent
Applying Lemma~\ref{lem:parabolic_estimates_matrix} to \(\delta f\), whose
initial datum is zero, yields
\[
\|\nabla\delta f(t)\|_{L^\infty}
\leq
C_T\|\psi\|_{\mathrm{BL}}\mathcal R_\tau\,t^{1/2},
\qquad 0<t\leq\tau.
\]
Therefore,
\begin{equation}
\label{eq:uniqueness-parabolic-estimate}
\int_0^\tau
\|\nabla\delta f(t)\|_{L^\infty}\,dt
\leq
\frac{2}{3}
C_T\|\psi\|_{\mathrm{BL}}
\tau^{3/2}\mathcal R_\tau.
\end{equation}
Finally, we estimate the difference of the velocity fields.
Subtracting the two velocity equations and writing the transport term
with respect to \(u_1\), we obtain
\[
\partial_t\delta u
+
(u_1\cdot\nabla)\delta u
=
-(\delta u\cdot\nabla)u_2
+
\mathcal A[\rho_1,u_1]
-
\mathcal A[\rho_2,u_2]
+
\gamma\nabla\delta f.
\]
Since \(\delta u(0)=0\), the \(L^\infty\)-estimate for the transport
equation gives, for every \(t\in[0,\tau]\),
\[
\begin{aligned}
\|\delta u(t)\|_{L^\infty}
&\leq
\int_0^t
\|\nabla u_2(s)\|_{L^\infty}
\|\delta u(s)\|_{L^\infty}\,ds
\\
&\quad+
\int_0^t
\left\|
\mathcal A[\rho_1(s),u_1(s)]
-
\mathcal A[\rho_2(s),u_2(s)]
\right\|_{L^\infty}\,ds+\gamma
\int_0^t
\|\nabla\delta f(s)\|_{L^\infty}\,ds.
\end{aligned}
\]
By Lemma~\ref{prop:alignment-estimates}, there exists a constant
\(C_A=C_A(M,m_0,\Phi)>0\) such that
$
\left\|
\mathcal A[\rho_1,u_1]
-
\mathcal A[\rho_2,u_2]
\right\|_{L^\infty}
\leq
C_A
\left[
d_{\mathrm{BL}}(\rho_1\,dx,\rho_2\,dx)
+
\|\delta u\|_{L^\infty}
\right].
$
Consequently,
\[
\|\delta u(t)\|_{L^\infty}
\leq
(M+C_A)
\int_0^t
\|\delta u(s)\|_{L^\infty}\,ds
+
C_A
\int_0^t
d_{\mathrm{BL}}(\rho_1(s)\,dx,\rho_2(s)\,dx)\,ds
+
\gamma
\int_0^t
\|\nabla\delta f(s)\|_{L^\infty}\,ds.
\]
Taking the supremum over \(t\in[0,\tau]\) and using
\eqref{eq:uniqueness-parabolic-estimate}, we obtain
\begin{equation}
\label{eq:uniqueness-velocity-estimate}
\mathcal U_\tau
\leq
(M+C_A)\tau\mathcal U_\tau
+
C_A\tau\mathcal R_\tau
+
\frac{2\gamma}{3}
C_T\|\psi\|_{\mathrm{BL}}
\tau^{3/2}\mathcal R_\tau.
\end{equation}
Inserting the density estimate
\eqref{eq:uniqueness-density-estimate}, into \eqref{eq:uniqueness-velocity-estimate}, we find
$
\mathcal U_\tau
\leq
\eta(\tau)\mathcal U_\tau,
$
where $\eta(\tau)
:=
(M+C_A)\tau
+
C_A m_0\tau^2e^{M\tau}
+
\frac{2\gamma}{3}
C_T\|\psi\|_{\mathrm{BL}}
m_0\tau^{5/2}e^{M\tau}$.\\

\noindent
Since \(C_T\) is fixed and finite and $\eta(\tau)\longrightarrow 0$ as $\tau\downarrow0$, there exists \(\tau_*\in(0,T]\), depending only on the uniform bounds of the two solutions and on the coefficients of the system,
such that $\eta(\tau_*)<1$. Therefore,
$
(1-\eta(\tau_*))\mathcal U_{\tau_*}\leq0.
$
Since \(\mathcal U_{\tau_*}\geq0\) and
\(1-\eta(\tau_*)>0\), it follows that $\mathcal U_{\tau_*}=0$.\\

\noindent
By the definition of \(\mathcal U_{\tau_*}\), this means that
$
u_1(t)=u_2(t)
\quad \text{in }L^\infty(\mathbb R^n;\mathbb R^n)
\quad\text{for every }t\in[0,\tau_*].
$
The density estimate
\eqref{eq:uniqueness-density-estimate} then yields $\mathcal R_{\tau_*}=0$, and hence
$
\rho_1(t,x)\,dx=\rho_2(t,x)\,dx
\qquad
\text{for every }t\in[0,\tau_*].
$
Since both measures are absolutely continuous with respect to the
Lebesgue measure, we conclude that
$
\rho_1(t)=\rho_2(t)
\quad
\text{in }L^1(\mathbb R^n)
\quad\text{for every }t\in[0,\tau_*].
$
Finally, \(\delta f=f_1-f_2\) satisfies
$
\partial_t\delta f
-
\operatorname{tr}
\bigl(\sigma\sigma^\top\nabla^2\delta f\bigr)
+
b\cdot\nabla\delta f
=
0
$
on \((0,\tau_*)\times\mathbb R^n\), with $\delta f(0)=0$. The uniqueness statement for the linear parabolic problem (see \cite{friedman2008partial}) therefore implies $f_1\equiv f_2$ on $\mathbb R^n \times [0,\tau_*]$. Consequently,
\[
(\rho_1,u_1,f_1)
=
(\rho_2,u_2,f_2)
\qquad\text{on }[0,\tau_*].
\]
The same argument can be restarted at \(t=\tau_*\). Indeed, the two
solutions have the same values at that time, and their norms remain
bounded by \(M\) on the whole interval \([0,T]\). Since all the constants entering the above estimates depend only 
on the common a priori bounds of the two solutions and not on the initial time, the argument can be restarted at 
$t=\tau_{*}$ with the same value of $\tau_{*}$. Repeating the argument a finite number of times, we conclude that
the solutions coincide on the entire interval $[0,T]$. This proves uniqueness in the class of Definition~\ref{def:strong-solution}.\\

\medskip
\noindent

The results obtained in Step 1 - Step 5 prove that \((\rho,u,f)\) satisfies all the six requirements of Definition~\ref{def:strong-solution}. This concludes the proof.
\end{proof}


\begin{remark}
\label{rem:regularity_class}
Strengthening \textup{(C1)} to $a,b\in C^{1+\alpha,(1+\alpha)/2}$ and assuming, in addition,
$\rho_0\in C^1_b(\mathbb R^n)$, $u_0\in C^2_b(\mathbb R^n)$, $f_0\in C^{3+\alpha}_b(\mathbb R^n)$,
$\Phi\in C^2_b(\mathbb R^n)$ and $\psi\in W^{2,\infty}(\mathbb R^n)$, the
same fixed-point argument can be carried out in the class
$\rho\in C^1(\mathbb R^n\times[0,T_*])$, $u,f\in C^{2,1}(\mathbb R^n\times[0,T_*])$. 
By the uniqueness part of that
theorem, the solution so obtained is the solution of Definition~\ref{def:strong-solution} issued
from the same data.
\end{remark}

\noindent
We conclude the section providing an a priori estimate for the velocity field in
the \(L^\infty\)-norm and establishing a continuation criterion for
maximal solutions. In particular, finite-time breakdown, if it occurs,
can only be caused by the loss of Lipschitz regularity of the velocity.\\

\begin{theorem}[Continuation criterion]
\label{thm:continuation-criterion}

Assume {\rm (H1)-(H3)} and {\rm (C1)-(C2)} on the time interval $[0,T]$, and set
$m_0:=\|\rho_0\|_{L^1}$. Let $(\rho,u,f)$ be the maximal solution in the sense of
Definition~\ref{def:strong-solution} issued from $(\rho_0,u_0,f_0)$, defined on
$[0,T_{\max})$ with $0<T_{\max}\le T$, whose existence and uniqueness are guaranteed by
Theorem~\ref{thm:local-wellposedness}. Then, for every \(\tau<T_{\max}\),
\begin{equation}
\label{eq:u-apriori-bound}
\sup_{0\leq t\leq \tau}\|u(t)\|_{L^\infty}
\leq
e^{2\|\phi\|_{L^\infty}m_0\tau}
\left[
\|u_0\|_{L^\infty}
+
\gamma C_T \tau\|\nabla f_0\|_{L^\infty}
+
\frac{2\gamma}{3}
C_T \|\psi\|_{L^\infty}m_0\tau ^{3/2}
\right],
\end{equation}
where \(C_T\) is the constant appearing in the parabolic estimates of
Lemma~\ref{lem:parabolic_estimates_matrix}. In particular, the \(L^\infty\)-norm of \(u\) cannot blow up in finite
time.

Moreover, if $T_{\max}<T$ and
$\sup_{0\leq t<T_{\max}}\|\nabla u(t)\|_{L^\infty}<+\infty$, then the solution can be extended to
an interval strictly larger than $[0,T_{\max})$. Equivalently, if $T_{\max}<T$, then
\begin{equation}
\label{eq:continuation_tesi}
\limsup_{t\uparrow T_{\max}}\|\nabla u(t)\|_{L^\infty}=+\infty .
\end{equation}
\end{theorem}

\begin{proof}\rm
We first establish the a priori bound
\eqref{eq:u-apriori-bound}.

\medskip
\noindent
By conservation of mass,
$\|\rho(t)\|_{L^1}=m_0$ for every $0\le t<T_{\max}$. Since $\rho(t,\cdot)\ge0$ almost
everywhere and $\|\Phi\|_{L^\infty(\mathbb R^n)}=\|\phi\|_{L^\infty(0,+\infty)}$, for every
$x\in\mathbb R^n$ the alignment term satisfies
\[
\bigl|\mathcal A[\rho,u](t,x)\bigr|
\leq
\int_{\mathbb R^n}
\Phi(x-y)\bigl(|u(t,y)|+|u(t,x)|\bigr)\rho(t,y)\,dy
\leq
2\|\phi\|_{L^\infty}m_0\|u(t)\|_{L^\infty}.
\]

Let $X(t;x)$ denote the characteristic flow generated by $u$, that is
$X'(t;x)=u\bigl(t,X(t;x)\bigr)$ with $X(0;x)=x$. For every $\tau<T_{\max}$ the triple $(\rho,u,f)$
is a solution on $[0,\tau]$ in the sense of Definition~\ref{def:strong-solution}, so that
$u\in L^\infty\bigl(0,\tau;W^{1,\infty}\bigr)$ and the flow is well defined on $[0,\tau]$.

Along the characteristics, the
velocity equation gives, for almost every $t$,
\[
\frac{d}{dt}u\bigl(t,X(t;x)\bigr)
=
\mathcal A[\rho,u]\bigl(t,X(t;x)\bigr)
+
\gamma\nabla f\bigl(t,X(t;x)\bigr).
\]

Integrating in time and taking the supremum with respect to $x$, we obtain
\[
\|u(t)\|_{L^\infty}
\leq
\|u_0\|_{L^\infty}
+
2\|\phi\|_{L^\infty}m_0\int_0^t\|u(s)\|_{L^\infty}\,ds
+
\gamma\int_0^t\|\nabla f(s)\|_{L^\infty}\,ds .
\]
By conservation of mass and estimate \eqref{eq:parabolic_gradient_estimate} with $M_\rho=m_0$, we have
\[
\int_0^\tau\|\nabla f(s)\|_{L^\infty}\,ds
\leq
C_T\tau\|\nabla f_0\|_{L^\infty}
+
\frac{2}{3}C_T\|\psi\|_{L^\infty}m_0\tau^{3/2},
\qquad \tau<T_{\max}.
\]
Since $t\mapsto\|u(t)\|_{L^\infty}$ is continuous and finite on $[0,\tau]$, Gr\"onwall's
inequality yields \eqref{eq:u-apriori-bound}.
We stress that the right-hand side depends only
on the data and on $\tau$, and not on $\|\nabla u\|_{L^\infty}$; in particular the
$L^\infty$-norm of $u$ remains bounded in $[0,T_{\max}) $.

\medskip
\noindent
We now turn to the second part. Let assume that $T_{\max}<T$ and that
\begin{equation}
\label{eq:K-finite}
K:=\sup_{0\leq t<T_{\max}}\|\nabla u(t)\|_{L^\infty}<+\infty .
\end{equation}
The characteristic representation of the density associated with the continuity equation reads
\[
\rho\bigl(t,X(t;x)\bigr)
=
\rho_0(x)\exp\left(-\int_0^t\operatorname{div}u\bigl(s,X(s;x)\bigr)\,ds\right),
\]
and, since $\|\operatorname{div}u(s)\|_{L^\infty}\leq n\|\nabla u(s)\|_{L^\infty}$, we obtain
\[
\|\rho(t)\|_{L^\infty}
\leq
\|\rho_0\|_{L^\infty}\exp\left(n\int_0^t\|\nabla u(s)\|_{L^\infty}\,ds\right)
\leq
\|\rho_0\|_{L^\infty}e^{nKT_{\max}}.
\]

Hence
\begin{equation}
\label{eq:rho-uniform-continuation}
\sup_{0\leq t<T_{\max}}\|\rho(t)\|_{L^\infty}<+\infty .
\end{equation}
Letting $\tau\uparrow T_{\max}$ in \eqref{eq:u-apriori-bound}, and combining the resulting bound
with the definition of $K$, we deduce
\begin{equation}
\label{eq:u-uniform-continuation}
\sup_{0\leq t<T_{\max}}\|u(t)\|_{W^{1,\infty}}<+\infty .
\end{equation}
As for the chemoattractant, the representation formula gives
$\|f(t)\|_{L^\infty}\le\|f_0\|_{L^\infty}+m_0\|\psi\|_{L^\infty}t$. 
Combining this with conservation of mass and estimates \eqref{eq:parabolic_gradient_estimate}-\eqref{eq:parabolic_hessian_estimate}, we obtain

\begin{equation}
\label{eq:f-uniform-continuation}
\sup_{0\le t<T_{\max}}\|f(t)\|_{W^{1,\infty}}<+\infty,
\qquad
\int_0^{T_{\max}}\|\nabla^2f(t)\|_{L^\infty}\,dt<+\infty .
\end{equation}

Consequently, all the norms entering the local existence theory remain bounded uniformly on
$[0,T_{\max})$.

\medskip
\noindent
Let $t_k\uparrow T_{\max}$ be any sequence. By \eqref{eq:rho-uniform-continuation},
\eqref{eq:u-uniform-continuation} and \eqref{eq:f-uniform-continuation}, the triples
$\bigl(\rho(t_k),u(t_k),f(t_k)\bigr)$ satisfy {\rm (H1)} with constants independent of $k$, and
$\|\rho(t_k)\|_{L^1}=m_0$. Moreover, both the existence and the uniqueness parts of
Theorem~\ref{thm:local-wellposedness} hold verbatim with initial time $t_k$ in place of $0$: the
constants in {\rm (C1)-(C2)} are uniform on $[0,T]$, and the constant $C_T$ of
Lemma~\ref{lem:parabolic_estimates_matrix} depends only on $n$, $T$ and on those constants. Since
the existence time is a nonincreasing function of the quantities listed in
Theorem~\ref{thm:local-wellposedness}, all of which are bounded uniformly in $k$, it can be chosen
equal to a single $\tau_0>0$, independent of $k$; recalling that $T_{\max}<T$, we may further
require $\tau_0\le T-T_{\max}$, so that $[t_k,t_k+\tau_0]\subseteq[0,T]$.

Fix $k$ so large that $t_k+\tau_0>T_{\max}$, and let $(\tilde\rho,\tilde u,\tilde f)$ be the
solution on $[t_k,t_k+\tau_0]$ issued from $\bigl(\rho(t_k),u(t_k),f(t_k)\bigr)$. By the
uniqueness statement of Theorem~\ref{thm:local-wellposedness}, applied on
$[t_k,T_{\max})$, the two solutions coincide there, so that the triple obtained by
concatenating $(\rho,u,f)$ on $[0,t_k]$ with $(\tilde\rho,\tilde u,\tilde f)$ on
$[t_k,t_k+\tau_0]$ is well defined and satisfies all the requirements of
Definition~\ref{def:strong-solution} on $[0,t_k+\tau_0]$; in particular
$\nabla^2f\in L^1\bigl(0,t_k+\tau_0;L^\infty(\mathbb R^n)\bigr)$, the two contributions being
separately integrable. 
This produces a solution on an interval strictly larger than
$[0,T_{\max})$, contradicting the maximality of $T_{\max}$. Hence
the supremum $K$ cannot be finite, proving \eqref{eq:continuation_tesi}.
\end{proof}

\section{Threshold conditions for finite-time breakdown}
\label{sec:multidimensional_thresholds}

The critical-threshold analysis follows the velocity-gradient matrix along the characteristics,
which requires differentiating the velocity equation pointwise in space. To this end, in this
section we assume
\[
\rho\in C^1(\mathbb R^n\times[0,\tau]),
\qquad
u,f\in C^{2,1}(\mathbb R^n\times[0,\tau]),
\]
for some $\tau\in(0,T]$, so that all the identities below hold in the classical sense. As
discussed in Remark~\ref{rem:regularity_class}, this class is reached by strengthening the
regularity of the coefficients and of the initial data, without affecting any of the estimates of
Theorem~\ref{thm:local-wellposedness}, and the solution so obtained coincides, by the uniqueness
part of that theorem, with the solution of Definition~\ref{def:strong-solution} issued from the
same data. Since Theorem~\ref{thm:continuation-criterion} shows that finite-time breakdown may
occur only through the blow-up of the velocity gradient, which is precisely the quantity estimated
below, the conditions derived in this section are conditions for the breakdown of that solution.

Let $(\rho,u,f)$ be a solution for the system \eqref{eq:model_intro} and, for a fixed $a\in \mathbb{R}^n$, we define $X(t;a)$ as the characteristic curve generated by the velocity field $u$, hence satisfying \eqref{eq:caratteristiche} with $\bar u=u$.
In the following, whenever a function $h=h(x,t)$ is evaluated along the characteristic issued from $a$, we use $h(t;a):=h\bigl(X(t;a),t\bigr)$.
We denote $M(x,t):=\nabla u(x,t)$ the velocity-gradient matrix and set $M(t):=M(X(t;a),t)$. Moreover, we define the divergence along characteristics by
\begin{equation}
d(t):=\operatorname{tr}M(t)=\nabla\cdot u(X(t;a),t).
\label{eq:_d}
\end{equation}

\noindent
We differentiate with respect to the spatial variable $x$ the second equation in \eqref{eq:model_intro}, getting
$
\nabla\big((u\cdot\nabla)u\big)
=
(u\cdot\nabla)M+M^2.
$
Hence
\begin{equation}
\partial_t M+(u\cdot\nabla)M+M^2
=
\nabla\mathcal A[\rho,u](x,t)+\gamma\nabla^2 f(x,t).
\label{eq:M_pde}
\end{equation}
Evaluating \eqref{eq:M_pde} along $X(t;a)$ gives
\begin{equation}
\frac{d}{dt}M(t)
=
\left(
\partial_t M+(u\cdot\nabla)M
\right)(X(t;a),t)
=
-M(t)^2
+
\nabla\mathcal A[\rho,u](X(t;a),t)
+
\gamma\nabla^2 f(X(t;a),t).
\label{eq:M_evolution}
\end{equation}
In order to compute the gradient of the alignment operator, we define $p(x,t):=
\int_{\mathbb R^n}\phi(|x-y|)\rho(y,t)\,dy$, and rewrite
\[
\mathcal A[\rho,u](x,t)
=
\int_{\mathbb R^n}
\phi(|x-y|)u(y,t)\rho(y,t)\,dy
-
p(x,t)u(x,t).
\]
Differentiating with respect to $x$, we get
\[
\nabla\mathcal A[\rho,u](x,t)
=
-p(x,t)\nabla u(x,t)
+
Q(x,t),
\]
where
\begin{equation}
Q(x,t)
:=
\int_{\mathbb R^n}
\big(u(y,t)-u(x,t)\big)\otimes\nabla\phi(|x-y|)\rho(y,t)\,dy.
\label{eq:Qal_def}
\end{equation}
Therefore, \eqref{eq:M_evolution} becomes
\begin{equation}
\frac{d}{dt}M(t)
=
-M(t)^2
-
p(t)M(t)
+
Q(t)
+
\gamma\nabla^2 f(X(t;a),t),
\label{eq:M_evolution_2}
\end{equation}
where $p(t):=p(X(t;a),t)$ and $Q(t):=Q(X(t;a),t)$. 
Taking the trace of \eqref{eq:M_evolution}, we obtain that the divergence along characteristics, defined in \eqref{eq:_d} satisfies
\begin{equation}
d'(t)
=
-\operatorname{tr}\big(M(t)^2\big)
-
p(t)d(t)
+
F_\gamma(t),
\label{eq:d_evolution_0}
\end{equation}
with 
\begin{equation}
\label{eq:Fgamma_def}
F_\gamma(t):=
\operatorname{tr}Q(t)
+
\gamma\Delta f(X(t;a),t).
\end{equation}

We decompose the trace as 
$
\operatorname{tr}\big(M(t)^2\big)
=
\frac1n d(t)^2+\mathcal S(t),
$
with
\begin{equation}
\mathcal S(t)
:=
\operatorname{tr}\big(M(t)^2\big)-\frac1n d(t)^2.
\label{eq:spectral_defect}
\end{equation}
Thus \eqref{eq:d_evolution_0} satisfies the scalar Riccati-type equation
\begin{equation}
d'(t)
=
-\frac1n d(t)^2
-
\mathcal S(t)
-
p(t)d(t)
+
F_\gamma(t).
\label{eq:d_evolution}
\end{equation}

\noindent
Before deriving the blow-up and non blow-up conditions, we establish uniform bounds on the alignment interaction term and the chemotactic forcing, which will play a crucial role in controlling the evolution of the velocity gradient along characteristics.

\begin{lemma}
\label{lem:bounds_p_Fgamma}
Let $(\rho,u,f)$ be a solution on $[0,T_*]$ such that 
$\sup_{t\in[0,T_*]}\|u(\cdot,t)\|_{L^{\infty}}\le R$ and $m_0:=\|\rho(\cdot,t)\|_{L^1}=\|\rho_0\|_{L^1}$. Then
\[
0\le p(x,t)\le B:=m_0 \|\phi\|_{L^\infty},\qquad\qquad\forall\; x\in \mathbb{R}^n,\;t\in [0,T_*].
\]
Moreover, the forcing term $F_\gamma$ in \eqref{eq:Fgamma_def} satisfies

\begin{equation}
|F_\gamma(t)|
\le 2n m_0 R \|\phi'\|_{L^\infty}+\gamma n C_{T_*}
\left(
\|\nabla f_0\|_{L^\infty}t^{-1/2}
+
m_0\|\nabla\psi\|_{L^\infty}t^{1/2}
\right)
\label{eq:Fgamma_bound_lemma}
\end{equation}
for any $0<t\le T_*$, where $C_{T_*}>0$ is the constant appearing in Lemma~\ref{lem:parabolic_estimates_matrix} and
depends only on $T_*$ and on the coefficients of the parabolic operator. 
\end{lemma}

\begin{proof}
Since $\phi\ge0$ and $\rho\ge0$,
$
0\le p(x,t)
\le
\|\phi\|_{L^\infty}
\int_{\mathbb R^n}\rho(y,t)\,dy
=
\|\phi\|_{L^\infty}m_0.
$
By the definition of $Q$ in (\ref{eq:Qal_def}), we get
\[
|\operatorname{tr}Q(x,t)|
\le
2n m_0 R \|\phi'\|_{L^\infty}.
\]
Applying Lemma \ref{lem:parabolic_estimates_matrix} with
$M_\rho=m_0$, we find
\[
\|\nabla^2f(\cdot,t)\|_{L^\infty}
\leq
C_{T_*}\|\nabla f_0\|_{L^\infty}t^{-1/2}
+
C_{T_*}m_0\|\nabla\psi\|_{L^\infty}t^{1/2},
\qquad 0<t\leq T_*.
\]
Consequently,
\[
\begin{aligned}
|F_\gamma(t)|
\leq {}&
2n m_0R\|\varphi'\|_{L^\infty}
+
\gamma n C_{T_*}
\left(
\|\nabla f_0\|_{L^\infty}t^{-1/2}
+
m_0\|\nabla\psi\|_{L^\infty}t^{1/2}
\right),
\end{aligned}
\]
for every $0<t\leq T_*$, proving \eqref{eq:Fgamma_bound_lemma}. 
\end{proof}

\noindent 
We now prove a sufficient condition for finite-time blow-up in the multidimensional setting. 
In the following we denote the positive and negative parts of $F_\gamma$ with
\[
F_\gamma^+(t):=\max\{F_\gamma(t),0\},
\qquad
F_\gamma^-(t):=\max\{-F_\gamma(t),0\}.
\]
Let $0<\delta<\tau<T$, and let $E_\delta^\gamma(\tau)$ and $I_\delta^\gamma$ be constants such that
\begin{equation}
\label{eq:nuovi:Edelta_Idelta}
F_\gamma^+(t)\le E_\delta^\gamma(\tau),
\quad \delta\le t\le\tau,
\qquad\qquad
\int_0^\delta F_\gamma^+(s)\,ds\le I_\delta^\gamma ,
\end{equation}
as long as the solution exists. By Lemma~\ref{lem:bounds_p_Fgamma}, together with the a priori
bound \eqref{eq:u-apriori-bound}, both constants may be chosen depending only on the
initial data and on the coefficients of the system, for fixed parameters $\delta$ and $\tau$. We set
\begin{equation}
\label{eq:nuovaKdelta}
K_\delta^\gamma(\tau)
:=
nB+\sqrt{n^2B^2+2nE_\delta^\gamma(\tau)},
\end{equation}
where $B$ is the constant of Lemma~\ref{lem:bounds_p_Fgamma}.

\begin{theorem}
\label{thm:multidim_breakdown_only}
Let $(\rho,u,f)$ be a maximal solution of the system on $[0,T_{\max})$. Assume that there exists $a\in\mathbb R^n$ such that
$\mathcal S(t)\ge0$ along the characteristic curve $X(\cdot;a)$, for every $t\in[0,T_{\max})$.
Let $0<\delta<\tau<T$, and let assume $E_\delta^\gamma(\tau)$, $I_\delta^\gamma\ge0$, $K_\delta^\gamma>0$ the constants in \eqref{eq:nuovi:Edelta_Idelta}, \eqref{eq:nuovaKdelta}. 
If the initial datum satisfies
\begin{equation}
-\nabla\cdot u_0(a)
\ge
e^{B\delta}
\left(
K_\delta^\gamma(\tau)
+
I_\delta^\gamma
\right),
\label{eq:supercritical_condition_nd_final}
\end{equation}
and
\begin{equation}
\delta+\frac{2n}{K_\delta^\gamma(\tau)}<\tau,
\label{eq:tau_condition_nd_final}
\end{equation}
then the maximal time of existence satisfies
\begin{equation}
T_{\max}
\le
\delta+\frac{2n}{K_\delta^\gamma(\tau)} .
\label{eq:Tmax_breakdown_nd_final}
\end{equation}
Consequently,
\[
\limsup_{t\uparrow T_{\max}}
\|\nabla u(\cdot,t)\|_{L^\infty}
=
+\infty .
\]
\end{theorem}

\begin{proof}\rm
Along the characteristic $X(\cdot;a)$, setting $r(t):=-d(t)$, the Riccati-type equation
\eqref{eq:d_evolution} reads
\[
r'(t)
=
\frac1n r(t)^2
+
\mathcal S(t)
-
p(t)r(t)
-
F_\gamma(t).
\]
Since $K_\delta^\gamma(\tau)>0$, assumption \eqref{eq:supercritical_condition_nd_final} gives in
particular $r(0)=-\nabla\cdot u_0(a)>0$.

\noindent
Suppose, by contradiction, that $r$ vanishes somewhere in $(0,\delta]$, and set
$t_*:=\inf\left\{t\in(0,\delta]\;:\;r(t)=0\right\}$. By continuity $r(t_*)=0$, and $r(t)>0$ for
every $t\in[0,t_*)$. Hence, using $\mathcal S\ge0$, $0\le p\le B$, $F_\gamma\le F_\gamma^+$ and
$\frac1nr^2\ge0$,
\[
r'(t)
\ge
-Br(t)-F_\gamma^+(t),
\qquad t\in[0,t_*),
\]
and the integrating factor $e^{Bt}$ yields
\[
r(t)\ge e^{-Bt}r(0)
-\int_0^t e^{-B(t-s)}F_\gamma^+(s)\,ds,
\qquad t\in[0,t_*].
\]
Since $t_*\le\delta$ and $e^{-B(t_*-s)}\le1$, letting $t\to t_*$ and using
\eqref{eq:supercritical_condition_nd_final} we obtain
\[
r(t_*)\ \ge\
e^{-B\delta}r(0)-I_\delta^\gamma
\ \ge\ K_\delta^\gamma(\tau)\ >\ 0,
\]
which contradicts $r(t_*)=0$. Hence $r>0$ on $[0,\delta]$, and the same estimate at $t=\delta$
gives $r(\delta)\ge K_\delta^\gamma(\tau)$.

\noindent
Assume now, by contradiction, that the solution exists on
$\bigl[\delta,\delta+2n/K^\gamma_\delta(\tau)\bigr]$; by \eqref{eq:tau_condition_nd_final} this
interval is contained in $[\delta,\tau]$, so that $F_\gamma^+\le E^\gamma_\delta(\tau)$ throughout.
We first show that $r(t)\ge K^\gamma_\delta(\tau)$ on that interval. Otherwise, let
$t_1:=\inf\bigl\{t>\delta\;:\;r(t)<K^\gamma_\delta(\tau)\bigr\}$; then
$r(t_1)=K^\gamma_\delta(\tau)$ and $r(t)<K^\gamma_\delta(\tau)$ for a sequence of times decreasing
to $t_1$, so that $r'(t_1)\le0$. On the other hand, since $K^\gamma_\delta(\tau)$ is the positive
root of $\frac1{2n}x^2-Bx-E^\gamma_\delta(\tau)$,
\[
r'(t_1)
\ge
\frac1n r(t_1)^2-Br(t_1)-E^\gamma_\delta(\tau)
\ge
\frac1{2n}r(t_1)^2>0,
\]
a contradiction. Consequently $r\ge K^\gamma_\delta(\tau)$ and, by the same computation,
$r'(t)\ge\frac1{2n}r(t)^2$ for every $t\in\bigl[\delta,\delta+2n/K^\gamma_\delta(\tau)\bigr]$.
Integrating $\bigl(1/r\bigr)'\le-\frac1{2n}$ from $\delta$ to $t$, we obtain
\[
\frac1{r(t)}
\le
\frac1{r(\delta)}-\frac{t-\delta}{2n}
\le
\frac1{K^\gamma_\delta(\tau)}-\frac{t-\delta}{2n},
\]
whose right-hand side vanishes at $t=\delta+2n/K^\gamma_\delta(\tau)$, which is impossible for a
finite positive $r$. This proves \eqref{eq:Tmax_breakdown_nd_final}. Finally, by
\eqref{eq:tau_condition_nd_final} we have $T_{\max}<\tau<T$, and
Theorem~\ref{thm:continuation-criterion} yields
$\limsup_{t\uparrow T_{\max}}\|\nabla u(\cdot,t)\|_{L^\infty}=+\infty$.

\end{proof}

\noindent
\begin{remark}\label{rmk_th_soglia} {\rm 1)} Note that if $n=1$, then $\mathcal S\equiv 0$. Therefore the assumption on the sign of $\mathcal S$  in the Theorem \ref{thm:multidim_breakdown_only} can be removed.\\
{\rm 2)} We emphasize that the threshold condition is formulated in terms of the initial velocity field $u_0$. However, the remaining initial data $\rho_0$ and $f_0$ enter the criterion indirectly through the quantities $p(\cdot)$ and $F_\gamma(\cdot)$. In particular, $\rho_0$ determines the total mass and thus the bound on the alignment term, while both $\rho_0$ and $f_0$ influence the forcing term $F_\gamma(\cdot)$ through the parabolic equation for $f$. 
\end{remark}

 While Theorem \ref{thm:multidim_breakdown_only} provides a sufficient condition for finite-time breakdown in the multidimensional setting, the following Theorem \ref{thm:1d_thresholds_gamma} establishes the complementary result in the one-dimensional case, where the scalar identity $M=\partial_x u$ makes the divergence control the entire velocity gradient.
In that case, the non blow-up conditions can be derived in an explicitly verifiable form, later on further investigated and validated with a numerical approach.

\begin{theorem}\label{thm:1d_thresholds_gamma}
Let $n=1$, and let $(\rho,u,f)$ be a maximal solution to \eqref{eq:model_intro} on
$[0,T_{\max})$, with $T_{\max}\le T$. 
For $a\in\mathbb R$ and $0\le t<T_{\max}$, 
let 
$
F_\gamma(t):=
F_\gamma(t;a)
$
as in \eqref{eq:Fgamma_def}. Let $0<\tau<T$, define
\[
Z_0:=\|(\partial_x u_0)^-\|_{L^\infty},
\qquad
(\partial_x u_0)^-:=\max\{0,-\partial_x u_0\},
\qquad
A_\tau^\gamma
:=
\sup_{a\in\mathbb R}
\int_0^{\min\{\tau,T_{\max}\}}
|F_\gamma(s;a)|\,ds .
\]
If
\begin{equation}
e^{B\tau}\tau\left(Z_0+A_\tau^\gamma\right)<1,
\label{eq:1d_nonblowup_condition}
\end{equation}
then $T_{\max}>\tau$ and
\begin{equation}
\sup_{0\le t\le \tau}
\|\partial_x u(\cdot,t)\|_{L^\infty}<+\infty.
\label{eq:bound_nabla_x_u}
\end{equation}
In particular, the strong solution can be continued at least beyond the time $\tau$.
\end{theorem}

\begin{proof}\rm (Theorem \ref{thm:1d_thresholds_gamma})
Set $\tau':=\min\{\tau,T_{\max}\}$. Since the solution is defined on $[0,T_{\max})$, all the
estimates below are performed on $[0,\tau')$, where every quantity involved is well defined

\noindent
Fix $a\in\mathbb R$. Since $n=1$, equation \eqref{eq:d_evolution} reduces to
$d'(t)=-d(t)^2-p(t)d(t)+F_\gamma(t)$. Define
$z(t):= z(t;a)=\max\{0,-d(t)\}=d^-(t)$ and $w(t):= w(t;a)=e^{-Bt}z(t)$. Since $p\ge0$ by
Lemma~\ref{lem:bounds_p_Fgamma}, in the sense of the upper Dini derivative we obtain
\[
D^+w(t)
\le
e^{B\tau}w(t)^2+F_\gamma^-(t),
\qquad 0\le t<\tau'.
\]
Integrating and applying Bihari-LaSalle's inequality yields
\begin{equation}
w(t)
\le
\frac{
z(0)+A_\tau^\gamma
}{
1-e^{B\tau}\big(z(0)+A_\tau^\gamma\big)t
},
\qquad 0\le t<\tau'.
\label{eq:w_a}
\end{equation}
Since $z(0)\le Z_0$ and $t<\tau'\le\tau$, condition \eqref{eq:1d_nonblowup_condition} bounds the
denominator in \eqref{eq:w_a} from below by $1-e^{B\tau}\tau\left(Z_0+A_\tau^\gamma\right)>0$,
uniformly with respect to $a\in\mathbb R$ and $t\in[0,\tau')$. Hence $w$, and therefore
$z=e^{Bt}w$, is bounded uniformly in $a$ and $t$. Recalling that the characteristic flow is a
bi-Lipschitz homeomorphism of $\mathbb R$ for every $t<T_{\max}$, this is equivalent to
\[
\sup_{0\le t<\tau'}
\|(\partial_xu(\cdot,t))^-\|_{L^\infty(\mathbb R)}<+\infty.
\]

\noindent
As for the positive part, whenever $d(t)>0$ we have $-d(t)^2-p(t)d(t)\le0$, and therefore
$d'(t)\le F_\gamma(t)\le F_\gamma^+(t)$. Integrating from the last time at which $d$
vanishes, or from $t=0$ if $d(\cdot)$ is positive on $[0,t]$, we obtain
\[
d(t)
\le
\max\{\partial_xu_0(a),0\}
+
\int_0^{\tau'}F_\gamma^+(s)\,ds
\le
\|\partial_xu_0\|_{L^\infty}+A_\tau^\gamma ,
\]
so that
\[
\sup_{0\le t<\tau'}
\|(\partial_xu(\cdot,t))^+\|_{L^\infty(\mathbb R)}<+\infty.
\]

Combining the estimates on the positive and negative parts leads to
\eqref{eq:bound_nabla_x_u}, which provides a uniform bound on the full velocity gradient. This is
precisely the quantity controlled by Theorem~\ref{thm:continuation-criterion}, and the strong
solution can therefore be continued beyond $\tau$. In particular, $T_{\max}>\tau$.

\end{proof}

\subsection{Analytical comparison between pressureless Euler alignment systems without and with chemotaxis }

With respect to classical pressureless Euler alignment systems, the presence of a chemotactic effect influencing the dynamics represents the main feature characterizing system \eqref{eq:model_intro}.
Unlike classical flocking models, the term $\gamma \nabla f$ is generated by the density through the parabolic equation and, via its integral representation, depends on the system history.
In this sense, chemotaxis here represents an endogenous memory mechanism in the alignment dynamics.

In this subsection we compare the chemotaxis-driven dynamics with the pure alignment model in the one-dimensional setting showing that, over finite times, the deviation from standard flocking behavior is controlled by the cumulative strength of the chemotactic field.

In one dimension, a sharp critical threshold has been established for the case without chemotaxis in \cite{carrillo2016critical}. Depending only on the initial configuration $(\rho_0, u_0)$, authors distinguish between \textit{subcritical} regions, ensuring global regularity, and \textit{supercritical} ones, leading to finite-time blow-up of the solution.  
The following Proposition \ref{prop:corrected_gradient_chemo} identifies the quantity that replaces the classical corrected-gradient invariant of the one-dimensional Euler-alignment dynamics studied in \cite{carrillo2016critical}, extending the result to the case of chemotaxis.

\begin{proposition}
\label{prop:corrected_gradient_chemo} Let \(n=1\), and let \((\rho,u,f)\) be the strong solution of system \eqref{eq:model_intro} on $[0,\tau]$. Define $e:=u_x+\phi*\rho$. Then along any characteristic \(t\mapsto X(t;a)\), satisfying $X^{\prime}(t;a)=u(X(t;a),t)$, $X(0;a)=a\in\mathbb{R}$, with $\rho_0(a)>0$, it holds:
\begin{equation}
\frac{e(X(t;a),t)}{\rho(X(t;a),t)}
=
\frac{u^\prime_0(a)+\phi*\rho_0(a)}{\rho_0(a)}
+
\gamma
\int_0^t
\frac{f_{xx}(X(s;a),s)}
{\rho(X(s;a),s)}
\,ds,\qquad\qquad \forall\;t \in [0,\tau].
\label{eq:e_rho}
\end{equation}
In particular, if \(\gamma=0\), the quantity \(e/\rho\) is preserved along characteristics.
\end{proposition}

\begin{proof}\rm
Differentiating the velocity equation with respect to \(x\), and setting \(v(x,t):=u_x(x,t)\), gives
\[
\partial_t v+u\,v_x+v^2
=
\partial_x\mathcal A[\rho,u]
+
\gamma f_{xx}.
\]
We now use the classical one-dimensional identity associated with the alignment operator. Since $\mathcal A[\rho,u]
=\phi*(\rho u)-u(\phi*\rho)$, we have
\[
\partial_x\mathcal A[\rho,u]
=
\partial_x\phi*(\rho u)
-
u_x(\phi*\rho)
-
u\,\partial_x(\phi*\rho).
\]
On the other hand, using the continuity equation, for $q:=(\phi*\rho)$, we obtain
\[
q_t=\phi*\rho_t
=
-\phi*(\rho u)_x
=
-\partial_x\phi*(\rho u)\quad \Longrightarrow\quad  q_t+u q_x
=
-\partial_x\phi*(\rho u)
+
u\,\partial_x(\phi*\rho).
\]
%
Adding the equations for \(v\) to the equation for \(q\), recalling that $e=v+q$, we get
\[
e_t+u e_x+v^2
=
-u_x q+\gamma f_{xx}.
\]
Since $v^2=u_x v$, this can be rewritten as
\[
e_t+u e_x
=
-u_x(v+q)+\gamma f_{xx}
=
-u_x e+\gamma f_{xx} \;\;\quad \Longrightarrow\;\;\quad e_t+(ue)_x=\gamma f_{xx}.\qquad\quad
\]

Evaluating along the characteristic $t\mapsto X(t;a)$, and using the characteristic representation
of the density associated with the continuity equation, we obtain
\begin{eqnarray*}
\frac{d}{dt}e(X(t;a),t)
&=&
-u_x(X(t;a),t)\,e(X(t;a),t)
+
\gamma f_{xx}(X(t;a),t),\\
\rho(X(t;a),t)&=& \rho_0(a)
\exp\left(
-\int_0^t u_x(X(s;a),s)\,ds
\right).
\end{eqnarray*}
In particular, since $\rho_0(a)>0$ and $\|u_x\|_{L^\infty(\mathbb R\times[0,\tau])}<+\infty$, the
second identity gives the lower bound
\begin{equation}
\label{eq:rho-lower-bound-char}
\rho(X(t;a),t)\ \ge\ \rho_0(a)\,e^{-\tau\|u_x\|_{L^\infty}}\ >\ 0,
\qquad t\in[0,\tau],
\end{equation}
so that $e/\rho$ is well defined along the characteristic. Differentiating it and using
$$\frac{d}{dt}\rho(X(t;a),t)=-u_x(X(t;a),t)\,\rho(X(t;a),t),$$ we compute
\[
\frac{d}{dt}
\left(
\frac{e}{\rho}
\right)(X(t;a),t)
=
\frac{e'\rho-e\rho'}{\rho^2}
=
\frac{
(-u_xe+\gamma f_{xx})\rho
-
e(-u_x\rho)
}
{\rho^2}
=
\gamma\,\frac{f_{xx}}{\rho}(X(t;a),t),
\]
where $e^\prime$ and $\rho^\prime$ denote the derivatives of $t\mapsto e(X(t;a),t)$ and
$t\mapsto\rho(X(t;a),t)$, respectively. By \eqref{eq:rho-lower-bound-char} and
\eqref{eq:parabolic_hessian_estimate}, the right-hand side is integrable on $(0,\tau)$, even
though $f_{xx}$ need not be bounded as $t\to0^+$. 
Integrating in time gives \eqref{eq:e_rho}. 
\end{proof}

We observe that, in the purely alignment-driven case (\(\gamma=0\)), 
the ratio \(e/\rho\) is conserved along characteristics. The quantity corresponds to function $\beta$ defined in \cite{carrillo2016critical}, proof of Proposition 2.1. When \(\gamma>0\), the ratio evolves according to a forcing term generated by the spatial curvature of the chemoattractant field. Consequently, the classical threshold mechanism becomes dynamically modified by the memory effects encoded in the parabolic equation. This implies that subcritical or supercritical behavior is therefore no longer determined solely by the initial configuration, but also by the cumulative contribution due to the past evolution of the density.

\subsection{The constant-coefficient, one-dimensional case}
\label{Constant-Coeff}

We now specialize the previous analysis to the one-dimensional setting with constant parabolic coefficients. Indeed, at present Euler-alignment systems with chemotaxis arise in modelling biological phenomena and typically account for constant diffusion and drift coefficients. Inspired by the literature \cite{menci2023microscopic, natalini2021mean, menci2026numericalMMS}, we assume
\begin{equation}
\label{eq:const_coeff_assumption}
\sigma(x,t)\equiv \sigma_0>0,
\qquad
b(x,t)\equiv b_0\in\mathbb R,
\qquad
f_0\equiv0.
\end{equation}

In this case, the fundamental solution of the parabolic equation is explicitly given by
\[
G(x,t;\xi,s)
=
\frac{1}{\sqrt{4\pi\sigma_0^2(t-s)}}
\exp\left(
-\frac{|x-\xi-b_0(t-s)|^2}{4\sigma_0^2(t-s)}
\right).
\]
Moreover, by the parabolic estimates established above, 
we get a sharp constant-coefficient estimate
\begin{equation}
\label{eq:hess_f_const_coeff}
\|\partial_{xx}f(\cdot,t)\|_{L^\infty}
\le
C_f\, t^{1/2},
\qquad
C_f:=\frac{2\,\|\psi'\|_{L^\infty}\,\|\rho_0\|_{L^1}}{\sigma_0\sqrt{\pi}},
\qquad 0<t\le T.
\end{equation}
In particular, $C_f$ depends only on the diffusion amplitude $\sigma_0$, on the total mass $\|\rho_0\|_{L^1}$, and on the Lipschitz constant $\|\psi'\|_{L^\infty}$ of the chemoattractant kernel, and is unaffected by the drift $b_0$. A sharper, variable-coefficient counterpart of $C_f$, retaining the explicit dependence on $b_0$ and $T$, can be deduced from Proposition~1 in \cite{menci2019global} through the corresponding estimate of $\partial_{xx}G$.

In the one-dimensional case, we recall that the forcing decomposes as
$F_\gamma(t) = Q(t) + \gamma\,\partial_{xx}f(X(t;a),t)$,
with $Q$ the alignment contribution of \eqref{eq:Qal_def}. Applying Lemma~\ref{lem:bounds_p_Fgamma} with $n=1$ and $f_0\equiv0$, and setting
\begin{equation}
\label{eq:CQ_def}
C_Q := 2\,\|\rho_0\|_{L^1}\,R\,\|\phi'\|_{L^\infty},
\qquad
R := \sup_{0\le t\le T}\|u(\cdot,t)\|_{L^\infty},
\end{equation}
the alignment bound $|Q(t)|\le C_Q$ combines with the parabolic estimate \eqref{eq:hess_f_const_coeff} to yield
\begin{equation}
\label{eq:Fgamma_1d_const}
|F_\gamma(t)|
\;\le\;
\underbrace{\,C_Q\,}_{\text{alignment}}
\;+\;
\underbrace{\,\gamma\, C_f\, t^{1/2}\,}_{\text{chemotaxis}},
\qquad \forall\; 0 < t \le T.
\end{equation}

\medskip
\noindent\textbf{Explicit blow-up and non blow-up thresholds}
Substituting \eqref{eq:Fgamma_1d_const} into Theorems~\ref{thm:multidim_breakdown_only} and~\ref{thm:1d_thresholds_gamma} turns the abstract criteria into explicit conditions on the initial slope. Since $n=1$,  $\mathcal S\equiv0$ (see Remark~\ref{rmk_th_soglia}), and $0\le p\le B$ with
$
B:=\|\rho_0\|_{L^1}\,\|\phi\|_{L^\infty}.
$

As the right-hand side of \eqref{eq:Fgamma_1d_const} is nondecreasing in $t$, admissible choices
of the constants entering Theorem~\ref{thm:multidim_breakdown_only} are, for $0<\delta<\tau\le T$,
\begin{eqnarray}
\label{eq_one_dim_est_blowup}
E_\delta^\gamma(\tau)
:= C_Q+\gamma C_f \tau^{1/2},
\qquad
K_\delta^\gamma(\tau)
= B+\sqrt{B^2+2\bigl(C_Q+\gamma C_f \tau^{1/2}\bigr)} .
\end{eqnarray}

Moreover, bounding the cumulative forcing over $[0,\delta]$ we get an explicit value for $I_\delta^\gamma$ in \eqref{eq:nuovi:Edelta_Idelta}. Indeed
\[
\int_0^\delta F_\gamma^+(s)\,ds
\le
\int_0^\delta \bigl(C_Q + \gamma C_f\, s^{1/2}\bigr)\,ds
=
C_Q\,\delta + \tfrac{2}{3}\gamma C_f\,\delta^{3/2}
=:
I_\delta^\gamma ,
\]

Substituting into \eqref{eq:supercritical_condition_nd_final}, a sufficient condition for finite-time breakdown along the characteristic issued from $a$ is
\begin{equation}
-u^\prime_0(a)
\ge
\mathcal{T}_{\delta,\tau}^\gamma,
\label{eq:supercritical_condition_nd_final_one}
\end{equation}
with the fully explicit threshold
\begin{equation}
\mathcal{T}_{\delta,\tau}^\gamma :=
e^{B\delta}
\left(
B+\sqrt{B^2+2\bigl(C_Q+\gamma C_f \tau^{1/2}\bigr)}
+
C_Q\,\delta+\frac{2}{3}\gamma C_f \,\delta^{3/2}\right).
\label{eq:TAU_gamma}
\end{equation}

If, in addition, the compatibility condition \eqref{eq:tau_condition_nd_final} holds, i.e.
$\delta+2/K_\delta^\gamma(\tau)<\tau$, then $T_{\max}\le\delta+2/K_\delta^\gamma(\tau)$ and, by
Theorem~\ref{thm:continuation-criterion},
$\limsup_{t\uparrow T_{\max}}\|\partial_xu(\cdot,t)\|_{L^\infty}=+\infty$.
Conversely, the continuation criterion of Theorem~\ref{thm:1d_thresholds_gamma} reduces to an explicit smallness condition on the negative part of the initial slope. Bounding the total forcing along characteristics,
\[
A_\tau^\gamma
\le
\int_0^\tau \bigl(C_Q + \gamma C_f\, s^{1/2}\bigr)\,ds
=
C_Q\,\tau + \frac{2}{3}\gamma C_f\,\tau^{3/2},
\]
uniformly in $a$, condition \eqref{eq:1d_nonblowup_condition} is implied by
\begin{equation}
\label{eq_one_dim_est_noblowup}
e^{B\tau}\,\tau
\left(
\|(\partial_xu_0)^-\|_{L^\infty}
+
C_Q\,\tau
+
\frac{2}{3}\gamma C_f\,\tau^{3/2}
\right)
<1 .
\end{equation}
We observe that both explicit thresholds \eqref{eq:TAU_gamma} and \eqref{eq_one_dim_est_noblowup} depend on the chemotactic strength $\gamma$ and reduce, for $\gamma=0$, to the corresponding thresholds of the pressureless Euler-alignment system in \cite{carrillo2016critical}.

\section{Numerical Results}
\label{sec:numerical_experiments}

In the following we numerically investigate pressureless Euler alignment systems with chemotaxis in the one-dimensional constant-coefficient regime. 
Specifically, the numerical approach focuses on a numerical study of the transition between finite-time breakdown and global regularity, and the role of chemotaxis.

We consider a computational domain $[-L,L]$, equipped with periodic boundary conditions, approximating systems \eqref{eq:model_intro} in the time interval $[0,T_{end}]$. Following the choice of previous works \cite{tadmor2014critical, albi2018pressureless}, we choose as communication kernel the function

\begin{equation}
\phi(r)
=
\frac{1}{(1+r^2)^{\alpha/2}},
\qquad\qquad
\alpha=0.5. 
\label{eq:phi_numeric}
\end{equation}
For the chemotactic interaction kernel, we consider
\begin{align}
\psi(x)
&=
\frac{1}{\sqrt{2\pi}\sigma_\psi}
\exp\left(
-\frac{x^2}{2\sigma_\psi^2}
\right).
\label{eq:psi_numeric}
\end{align}
which models a localized smooth chemoattractant production mechanism. 

We compute the $ L^\infty$-norms for the chemoattractant kernel and its derivative, required in the following computations, getting

\begin{equation}
\label{eq:psi_norms}
\|\psi\|_{L^\infty}=\frac{1}{\sqrt{2\pi}\,\sigma_\psi},
\qquad
\|\psi'\|_{L^\infty}=\frac{e^{-1/2}}{\sqrt{2\pi}\,\sigma_\psi^{2}} .
\end{equation}

In all the numerical simulations, we consider as initial conditions 
\begin{align}
\rho_0(x)
&=
\frac{1}{\sigma_\rho\sqrt{2\pi}}
\exp\!\left(-\frac{x^2}{2\sigma_\rho^2}\right),
\qquad \sigma_\rho = 0.3,\quad \|\rho_0\|_{L^1}=1,
\label{eq:rho0_test1D}
\\
u_0(x) &= -A\sin\!\left(\frac{\pi x}{L}\right),
\label{eq:u0_test1D}
\\
f_0(x) &= 0.
\label{eq:f0_test1D}
\end{align}

We observe that the velocity initial profile is $C^\infty$ and $2L$-periodic, and it is compressive near the origin, with 
$\partial_x u_0(0)=-A\pi/L$, $A>0$.
%
%
In line with numerical investigations performed in \cite{menci2026numericalMMS}, we here set $b_0=0$.
Under these assumptions the two parabolic estimates on $\partial_x f$ and  $\partial_{xx}f$ can be rewritten as
\begin{equation}
\label{eq:two_parabolic_const}
\|\partial_x f(\cdot,t)\|_{L^\infty}\le \tilde C_f\,t^{1/2},
\qquad
\|\partial_{xx}f(\cdot,t)\|_{L^\infty}\le C_f\,t^{1/2},
\qquad 0<t\le T,
\end{equation}
with the two distinct constants
\begin{equation}
\label{eq:Cf_two}
\tilde C_f:=\frac{2\,\|\psi\|_{L^\infty}\,\|\rho_0\|_{L^1}}{\sigma_0\sqrt{\pi}},
\qquad
C_f:=\frac{2\,\|\psi'\|_{L^\infty}\,\|\rho_0\|_{L^1}}{\sigma_0\sqrt{\pi}}
=\frac{e^{-1/2}}{\sigma_\psi}\,\tilde C_f .
\end{equation}

Since $\phi(0)=1$, we have $B=\|\rho_0\|_{L^1}\|\phi\|_{L^\infty}=1$. In the generic regime
\eqref{eq:Fgamma_1d_const} reads
\begin{equation}
\label{eq:Fgamma_numeric}
|F_\gamma(t)|\le C_Q+\gamma C_f\,t^{1/2},
\qquad
C_Q=2\,\|\rho_0\|_{L^1}\,R\,\|\phi'\|_{L^\infty},
\qquad
R=\sup_{0\le t\le\tau}\|u(\cdot,t)\|_{L^\infty}.
\end{equation}

\medskip
We observe that evaluating $C_Q$ requires an estimate of the velocity amplitude $R$, we here compute considering a maximum principle approach applied to equation \eqref{eq:model_intro}$_2$. 
Along the characteristics $\dot x=u$, one has $\tfrac{d}{dt}u=\mathcal A[\rho,u]+\gamma\,\partial_x f$. At a point where $u(\cdot,t)$ attains its spatial maximum, the alignment operator is non-positive, $\mathcal A[\rho,u]=\int\phi(|x-y|)\bigl(u(y)-u(x)\bigr)\rho(y)\,dy\le0$, since $\phi,\rho\ge0$ and $u(y)\le u(x)$. With symmetric argumentation, it is non-negative at the minimum.
Hence we get
\[
\frac{d}{dt}\max_x u(\cdot,t)\le \gamma\,\|\partial_x f(\cdot,t)\|_{L^\infty},
\qquad
\frac{d}{dt}\min_x u(\cdot,t)\ge -\gamma\,\|\partial_x f(\cdot,t)\|_{L^\infty}.
\]
In other words, the alignment term can only contract the velocity range, and the mechanism responsible for enlarging it is the chemotactic drift.
Integrating in time and inserting the gradient estimate \eqref{eq:two_parabolic_const},
\begin{equation}
\label{eq:R_estimate}
\|u(\cdot,t)\|_{L^\infty}
\le \|u_0\|_{L^\infty}+\gamma\int_0^t\|\partial_x f(\cdot,s)\|_{L^\infty}\,ds
\le \|u_0\|_{L^\infty}+\tfrac{2}{3}\,\gamma\,\tilde C_f\,t^{3/2}.
\end{equation}
We observe that the right-hand side is increasing in $t$ and, since in the supercritical regime breakdown occurs at some $T_{\max}\le\tau$, its value at $t=\tau$ is chosen as the a priori upper bound for $R$. In conclusion, we obtain
\begin{equation}
\label{eq:R_apriori}
R\le \|u_0\|_{L^\infty}+\tfrac{2}{3}\,\gamma\,\tilde C_f\,\tau^{3/2},
\qquad
C_Q\le 2\,\|\rho_0\|_{L^1}\,\|\phi'\|_{L^\infty}\Bigl(\|u_0\|_{L^\infty}+\tfrac{2}{3}\,\gamma\,\tilde C_f\,\tau^{3/2}\Bigr).
\end{equation}
For our initial condition \eqref{eq:u0_test1D}, $\|u_0\|_{L^\infty}=A$, so \eqref{eq:R_apriori} makes $C_Q$ explicitly computable from the data. 
With $B=1$ and the bound \eqref{eq:Fgamma_numeric}, the admissible choices
\eqref{eq_one_dim_est_blowup} and the threshold \eqref{eq:TAU_gamma} specialize, in the numerical
setting, to
\begin{align}
K^\gamma_\delta(\tau)
&=
1+\sqrt{\,1+2\bigl(C_Q+\gamma C_f\,\tau^{1/2}\bigr)}\,,
\label{eq:Kgamma_b0}
\\[4pt]
\mathcal{T}_{\delta,\tau}^\gamma
&=
e^{\delta}
\Bigl(
1+\sqrt{\,1+2\bigl(C_Q+\gamma C_f\,\tau^{1/2}\bigr)}
+C_Q\,\delta+\tfrac{2}{3}\gamma C_f\,\delta^{3/2}
\Bigr).
\label{eq:Tthreshold_b0}
\end{align}

\medskip
Based on the sufficient condition \eqref{eq:supercritical_condition_nd_final_one}, we define the \textit{threshold ratio}
\begin{equation}
\mathcal{R} := \frac{-\inf_x u_0'(x)}{\mathcal{T}_{\delta,\tau}^\gamma}.
\label{eq:R_ratio}
\end{equation}

We highlight that this ratio does not define a sharp critical threshold, since our
criterions are only sufficient.
For this reason, we adapt the definition introduced in \cite{carrillo2016critical} to our case. In the following, we refer to the following behaviors in this sense:
\begin{itemize}
\item \emph{Subcritical case} ($\gamma = \gamma_{\rm sub}$, $\mathcal{R}<1$): the sufficient condition ~\eqref{eq_one_dim_est_noblowup} is satisfied, and the gradient $\|\partial_x u(\cdot,t)\|_{L^\infty}$ remains bounded over $[0,T_{\rm end}]$.
\item \emph{Supercritical case} ($\gamma = \gamma_{\rm sup}$, $\mathcal{R}>1$): the blow-up condition~\eqref{eq:supercritical_condition_nd_final_one} holds, and the simulation detects finite-time blow-up at a numerical time $T_{\rm num}$ (defined in the following \eqref{eq:Tnum}) consistent with the theoretical estimate \eqref{eq:Tmax_breakdown_nd_final}.
\end{itemize}

In our setting, we get that
\begin{equation}
\mathcal{R}
=
\frac{-\partial_x u_0(0)}{\mathcal{T}_{\delta,\tau}^\gamma}
=
\frac{A\pi/L}{\mathcal{T}_{\delta,\tau}^\gamma}.
\label{eq:R}
\end{equation}

\paragraph{The numerical protocol} 
First, we evaluate the explicit threshold quantities  $\mathcal{T}_{\delta,\tau}^\gamma$ together with the associated maximal theoretical time 
\begin{equation}
T_{\mathrm{UB}}(\tau):=
\delta + \frac{2}{K^\gamma_\delta(\tau)}.
\label{eq:TUB_b0}
\end{equation}
In order to numerically exploit the behavior of the solutions we monitor the quantities
\begin{equation}
\label{eq:M_m}
M_u(t):=\|\partial_xu(\cdot,t)\|_{L^\infty}, 
\qquad 
m_u(t):=\min_x \partial_xu(x,t).
\end{equation}
The numerical breakdown time is defined as
\begin{equation}
T_{\rm num}
:=
\inf\left\{
t>0:
\|\partial_xu(\cdot,t)\|_{L^\infty}\ge M_{\rm tol}
\right\},
\label{eq:Tnum}
\end{equation}
where $M_{\rm tol}$ is a fixed tolerance. Finally, we compare $T_{\rm num}$ with the theoretical estimate $T_{UB}$. In the supercritical regime, $m_u(t)\to -\infty$ and $M_u(t)$ rapidly growth as $t$ approaches the detected breakdown time, whereas in the subcritical regime $M_u(t)$ should remain bounded over the prescribed time interval.

\paragraph{Simulation setting}

For the spatial domain, we fix $L=1$ and a discretization step $\Delta x= 0.005$, corresponding to $N=400$ discretization points.
The simulation runs over the time interval $[0,T]=[0,1.5]$, with $\Delta t$ chosen according to the
CFL condition required by the scheme adopted (see the following \eqref{eq:CFL}), choosing $C_{\rm CFL}=0.4$.
We fix the parabolic diffusion coefficient $\sigma_0=0.1$, and the chemotactic coupling $\gamma=0.5$ in Test 1, whereas $\gamma$ varies in Test 2. 

For the communication kernel \eqref{eq:phi_numeric} with $\alpha=0.5$, 
$\|\phi'\|_{L^\infty}=\max_{r\ge0}\tfrac{r}{2}(1+r^2)^{-5/4}\approx0.2156$.
Analogously, from \eqref{eq:psi_norms} choosing $\sigma_\psi=0.3$, we get $\|\psi\|_{L^\infty}\approx1.330$ and $\|\psi'\|_{L^\infty}\approx2.689$.

For each value of $\gamma$, the alignment constant $C_Q$ is then obtained from \eqref{eq:R_apriori}, the quantities $\mathcal{T}_{\delta,\tau}^\gamma$,
$T_{\rm UB}$, and the threshold ratio $\mathcal{R}$
are computed from \eqref{eq:Tthreshold_b0}, \eqref{eq:TUB_b0}, \eqref{eq:R_ratio} respectively,
with $\delta=0.05$, $\tau=0.5$ and different values of $A$ reported in the following tables.

The relevant constant entering the evaluation of  $F_\gamma$, $K^\gamma_\delta$ and $\mathcal{T}_{\delta,\tau}^\gamma$ is $C_f\approx30.34$, whereas $\tilde C_f\approx15.00$ is required in the evaluation of the bound \eqref{eq:R_apriori}.
Both values are computed from \eqref{eq:Cf_two}, using $\|\rho_0\|_{L^1}=1$.

\subsection{A direct finite-difference scheme}
\label{subsec:scheme2}

\noindent
In the following we give details on the finite-difference numerical scheme implemented to discretize system \eqref{eq:model_intro}, in our numerical setting.
We consider the Crank-Nicolson scheme to approximate the solution $f$ at any grid point $i=1, ..., N$:
\begin{equation}
\frac{f_i^{n+1}-f_i^n}{\Delta t}
=
\frac{\sigma_0^2}{2}
\!\left(D^2_x f_i^{n+1}+D^2_x f_i^n\right)
-
\frac{b_0}{2}
\!\left(D_x f_i^{n+1}+D_x f_i^n\right)
+
(\psi*\rho^n)_i,
\label{eq:CN_f}
\end{equation}
where $D^2_x f_i = (f_{i+1}-2f_i+f_{i-1})/\Delta x^2$ and $D_x f_i=(f_{i+1}-f_{i-1})/(2\Delta x)$ are the standard second- and first-order central-difference operators.  Rearranging~\eqref{eq:CN_f}, the implicit step amounts to solving the linear system
\begin{equation}
\left(I - \tfrac{\Delta t}{2}\mathcal{L}_h\right)f^{n+1}
=
\left(I + \tfrac{\Delta t}{2}\mathcal{L}_h\right)f^n
+
\Delta t\,(\psi*\rho^n),
\label{eq:CN_system}
\end{equation}
Assuming $b_0=0$, the discrete operator $\mathcal{L}_h$
reduces to the pure second-order finite-difference Laplacian
$\mathcal{L}_h = \sigma_0^2 D^2_x$, and the first row of the circulant
matrix $(I-\tfrac{\Delta t}{2}\mathcal{L}_h)$ is
\[
\bigl[1+r,\;-r/2,\;0,\;\ldots,\;0,\;-r/2\bigr], \qquad r=\frac{\sigma_0^2\,\Delta t}{\Delta x^2}.
\]
We observe that, with periodic boundary conditions and constant coefficients, both matrices
$(I\pm\tfrac{\Delta t}{2}\mathcal{L}_h)$ are circulant.  Consequently, system~\eqref{eq:CN_system} is solved exactly by one forward and one inverse Fast Fourier Transform (FFT):
\begin{equation}
f^{n+1}
=
\mathcal{F}^{-1}\!\left(\frac{\mathcal{F}(rhs^n)}{\hat\lambda_{\mathrm{lhs}}}\right),
\label{eq:CN_FFT}
\end{equation}
where $\hat\lambda_{\mathrm{lhs}}$ is the vector of eigenvalues of $(I-\tfrac{\Delta t}{2}\mathcal{L}_h)$, given by the FFT of its first row, and $rhs^n$ denotes the right-hand side in~\eqref{eq:CN_system}. With this numerical choice we avoid any direct matrix inversion.

The conservation law for $\rho$ and the velocity equation for $u$ are advanced with explicit upwind stepping, using the gradient $\partial_x f^n$ (evaluated at the previous step) as the chemotactic forcing:
\begin{align}
u_i^{n+1}
&=
u_i^n
-
\Delta t\,u_i^n D_xu_i^n
+
\Delta t\,\mathcal A_i^n
+
\Delta t\,\gamma \partial_x f_i^n,
\label{eq:u_scheme2}
\\
\rho_i^{n+1}
&=
\rho_i^n
-
\Delta t\,D_x(\rho^n u^n)_i.
\label{eq:rho_scheme2}
\end{align}

The discrete alignment term $\mathcal{A}_i^n$ is the Cucker-Smale interaction
term evaluated at grid point $i=1,...,N$ and time $t_n$:
\begin{equation}
\mathcal{A}_i^n
\;=\;
\Delta x\sum_{j=1}^{N}\phi\!\left(\lvert x_i - x_j\rvert\right)\rho_j^n(u_j^n-u_i^n)
\;=\;
\bigl(\phi*(\rho^n u^n)\bigr)_i \;-\; u_i^n\,\bigl(\phi*\rho^n\bigr)_i,
\label{eq:alignment_discrete}
\end{equation}
with $\phi$ defined in \eqref{eq:phi_numeric}.
Both convolutions in~\eqref{eq:alignment_discrete} are computed in $\mathcal{O}(N\log N)$ via FFT.  
The hyperbolic updates~\eqref{eq:u_scheme2}-\eqref{eq:rho_scheme2} are subject to the standard hyperbolic CFL condition
\begin{equation}
\Delta t \le C_{\mathrm{CFL}}\,\frac{\Delta x}{\max_i |u_i^n|}.
\label{eq:CFL}
\end{equation}

At each step we monitor the quantities in \eqref{eq:M_m}, numerically discretized as:
\begin{equation}
M_u(t_n) := \|\partial_x u(\cdot,t_n)\|_{L^\infty}
\approx
\max_{i}\Bigl|\frac{u_{i+1}^n - u_{i-1}^n}{2\Delta x}\Bigr|,
\qquad
m_u(t_n) := \min_{i}\frac{u_{i+1}^n - u_{i-1}^n}{2\Delta x}.
\label{eq:Mu_monitor}
\end{equation}
In the following Test 1 and Test 2, the numerical breakdown time as defined in \eqref{eq:Tnum}
is evaluated fixing $M_{\rm tol}=100$.

\subsection{Test 1: behavior of the solution in different regimes}

In this test, performed at fixed $\gamma=0.5$, we simulate both a subcritical and a supercritical scenario.
Numerical results are collected in Table \ref{tab:levels_test1}. We remark that the theoretical bound $T_{\rm UB}$ is meaningful only in the supercritical case $\mathcal{R}>1$. Indeed, in the subcritical case, the blow-up criterion is not met and no breakdown is expected within $T_{\rm end}$.

\begin{table}[H]
\centering
\caption{Numerical Results of Test1.}
\label{tab:levels_test1}

\begin{tabular}{lcccccc}
\toprule
Regime & $A$  & $\mathcal{R}$ & $-\partial_x u_0(0)=A\pi/L$ & $\mathcal{T}_{\delta,\tau}^\gamma$ & $T_{\rm num}$ & $T_{\rm UB}$ \\
\midrule
Subcritical    & $A_{\rm sub}=0.10$ & $0.05$ & $0.3247$ & $6.37$ & $>1.5$  & - \\
Supercritical &  $A_{\rm sup}=3.30$ & $1.54$ & $10.388$ & $6.73$ & $0.0950$  & $0.3736$ \\
\bottomrule
\end{tabular}
\end{table}

In the subcritical case the initial compression $-\partial_x u_0(0)=0.32$ is only about $5\%$ of the threshold $\mathcal{T}_{\delta,\tau}^\gamma\approx6.37$, so that $\mathcal{R}\approx0.051$. The initial velocity gradient is extremely small and the alignment interaction dominates: no blow-up is detected within $T_{\rm end}=1.5$, consistently with $\mathcal{R}\ll1$.

In the supercritical case the initial compression $-\partial_x u_0(0)=10.39$ exceeds the threshold $\mathcal{T}_{\delta,\tau}^\gamma\approx6.73$, giving $\mathcal{R}\approx1.54>1$. The inward compression overwhelms both the alignment and chemotactic effects, and the velocity gradient $M_u(t)$ grows and crosses $M_{\rm tol}=100$ at $T_{\rm num}=0.0950$. The consistency check $T_{\rm num}=0.0950\le T_{\rm UB}=0.3736$ is satisfied, numerically confirming the validity of the theoretical upper bound.

\begin{figure}[H]
\centering
\subfigure[Density $\rho$ ($\mathcal{R}\approx0.05$)]{%
  \includegraphics[width=0.46\textwidth]{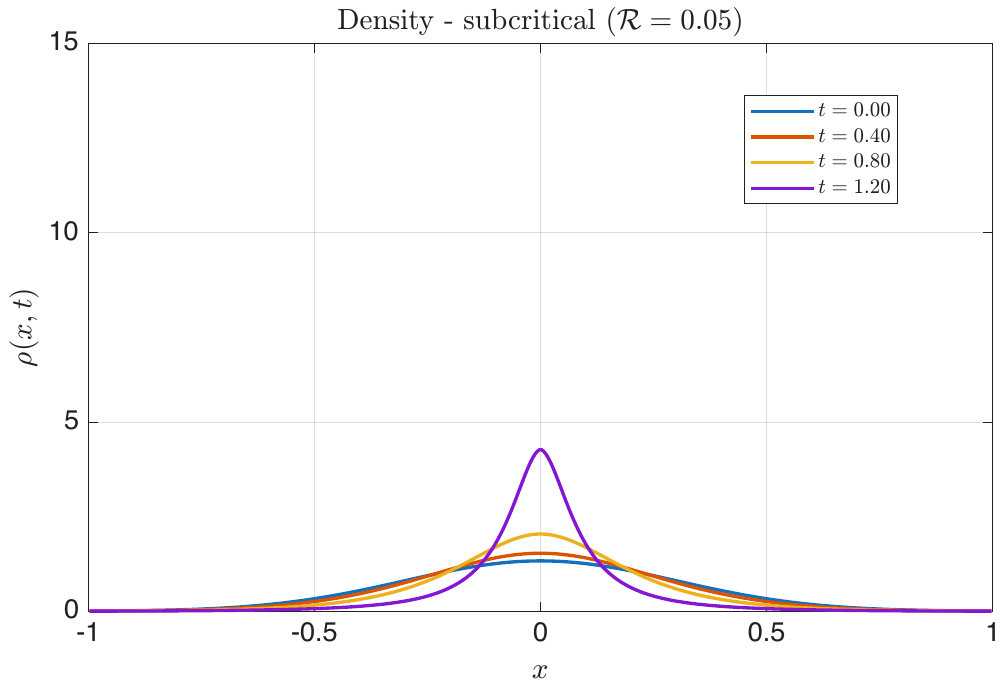}}\hfill
\subfigure[Density $\rho$ ($\mathcal{R}\approx1.54$)]{%
  \includegraphics[width=0.46\textwidth]{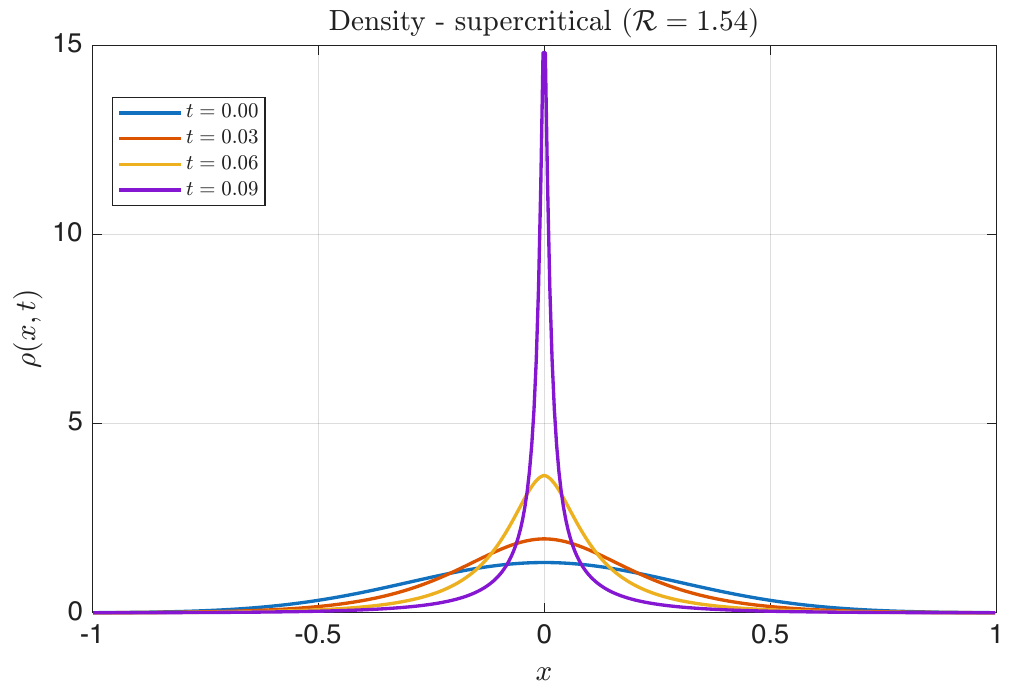}}\\[6pt]
\subfigure[Velocity $u$ ($\mathcal{R}\approx0.05$)]{%
  \includegraphics[width=0.46\textwidth]{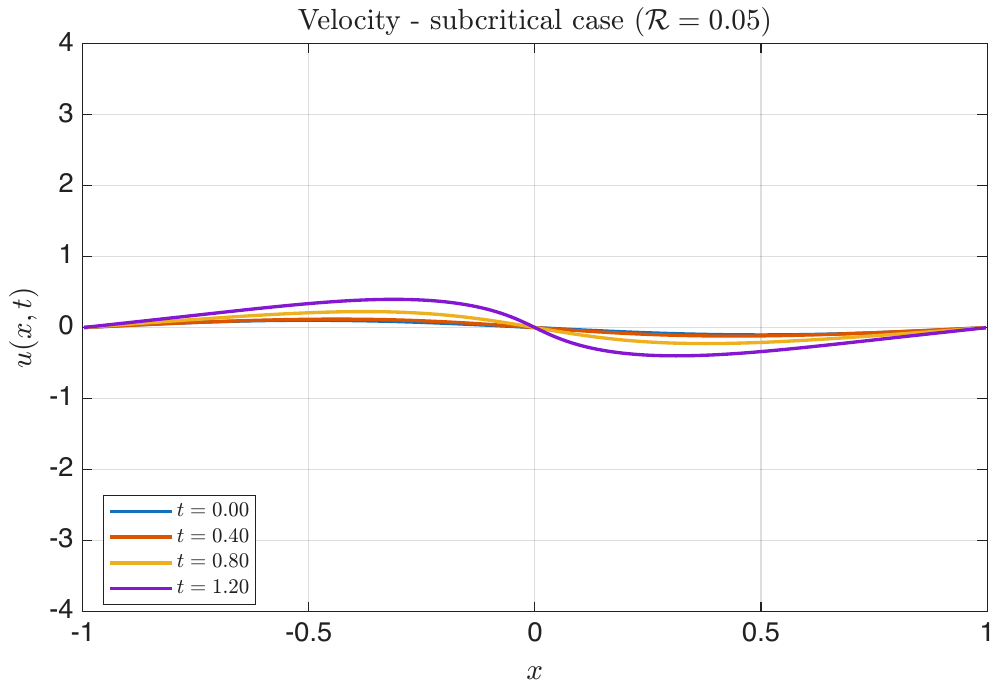}}\hfill
\subfigure[Velocity $u$ ($\mathcal{R}\approx1.54$)]{%
  \includegraphics[width=0.46\textwidth]{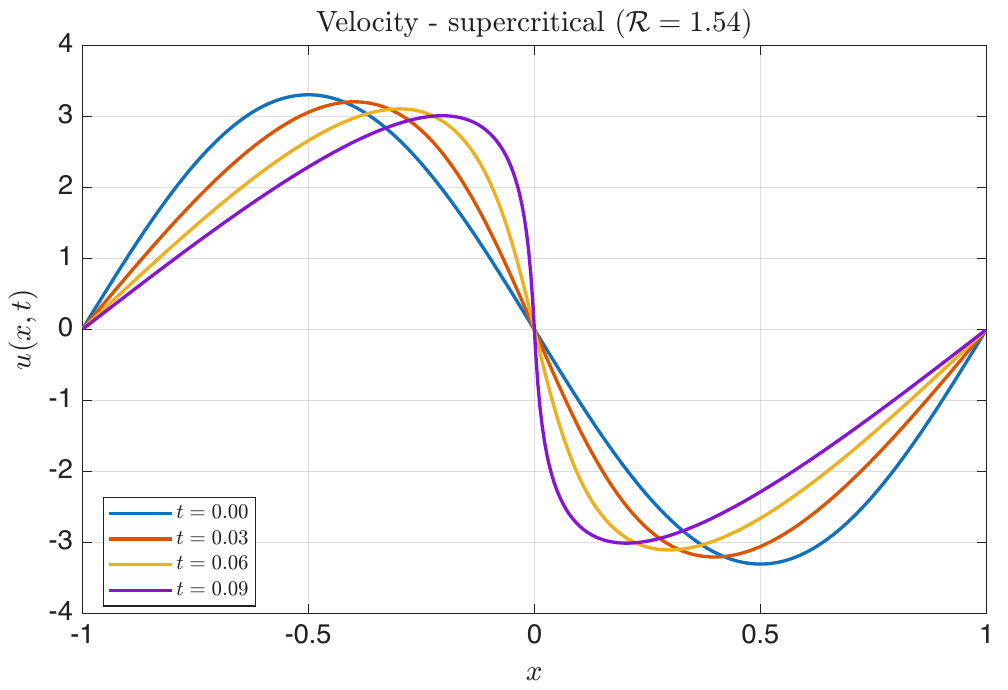}}
\caption{%
Comparison of density and velocity snapshots between the subcritical
($\mathcal{R}\approx0.05$) and supercritical ($\mathcal{R}\approx1.54$) regimes.
(a)-(b): plot of density $\rho(x,t)$ at different times $t$; (c)-(d): plot of velocity $u(x,t)$ at different times $t$.
}
\label{fig:panel_sub_sup}
\end{figure}

Figure \ref{fig:panel_sub_sup}  illustrates the qualitative difference between the two regimes.
In the subcritical case (left column), the density profile remains smooth over a long time interval and spreads only mildly, while the velocity gradient grows slowly.
In the supercritical  case (right column),  the density concentrates sharply near the origin and the
velocity field develops a steep front within a very short time, with blow-up at $T_{\rm num}=0.0950$.

\paragraph{Sensitivity of detected blow-up.}
We conclude the analysis of Test~1 by assessing the robustness of the detected blow-up time $T_{\rm num}$ with respect to the discretization and detection parameters, in the supercritical scenario. To this end we vary the number of discretization points $N$, the time step, through the $C_{\rm CFL}$ parameter, and the blow-up tolerance $M_{\rm tol}$. 
The analysis is performed varying one parameter at a time, holding the other two at their reference values, i.e. $N=400$, $C_{\rm CFL}=0.4$, $M_{\rm tol}=100$. The results obtained are collected in Table \ref{tab:sensitivity}.

\begin{table}[H]
\centering
\caption{Sensitivity of the detected blow-up time $T_{\rm num}$ in the supercritical regime of Test~1. Each block varies a single parameter, with the other two fixed at the baseline $N=400$, $C_{\rm CFL}=0.4$, $M_{\rm tol}=100$.}
\label{tab:sensitivity}
\small
\begin{tabular}{cc}
\toprule
\multicolumn{2}{c}{(a) mesh size}\\
\midrule
$N$ & $T_{\rm num}$ \\
\midrule
$100$  & $0.1185$ \\
$200$  & $0.1004$ \\
$400$  & $0.0950$ \\
$800$  & $0.0925$ \\
$1600$ & $0.0916$ \\
\bottomrule
\end{tabular}
\begin{tabular}{cc}
\toprule
\multicolumn{2}{c}{(b) time step}\\
\midrule
$C_{\rm CFL}$ & $T_{\rm num}$ \\
\midrule
$0.80$ & $0.0968$ \\
$0.40$ & $0.0950$ \\
$0.20$ & $0.0941$ \\
$0.10$ & $0.0938$ \\
$0.05$ & $0.0936$ \\
\bottomrule
\end{tabular}
\begin{tabular}{cc}
\toprule
\multicolumn{2}{c}{(c) tolerance}\\
\midrule
$M_{\rm tol}$ & $T_{\rm num}$ \\
\midrule
$50$   & $0.0829$ \\
$100$  & $0.0950$ \\
$200$  & $0.1034$ \\
$500$  & $0.1240$ \\
$1000$ & $>1.5$ \\
\bottomrule
\end{tabular}

\end{table}

Under mesh refinement (Table~\ref{tab:sensitivity}a) the detected time decreases monotonically and converges. The first-order numerical diffusion of the upwind scheme, of size $\sim|u|\,\Delta x/2$, mildly delays the steepening on coarse grids, and its effect vanishes as $N$ increases.
The sensitivity to the time step (Table~\ref{tab:sensitivity}b) is even weaker. Indeed we observe that reducing $\Delta t$ by a factor $16$ changes $T_{\rm num}$ by around $3\%$, consistent with the spatial error dominating the temporal one.

The dependence on the detection tolerance (Table~\ref{tab:sensitivity}c) is more pronounced and the numerical analysis gives interesting insights.
The time $T_{\rm num}$ increases with $M_{\rm tol}$ ($0.0829$ at $M_{\rm tol}=50$ up to $0.1240$ at $M_{\rm tol}=500$), since a larger gradient threshold is crossed slightly later along the rapid growth of $M_u$. Moreover, on a fixed grid the discrete gradient cannot grow unboundedly: numerical diffusion caps it at a resolution-dependent ceiling $\max_{[0,T_{\rm end}]}M_u\approx C/\Delta x$, which scales linearly with $N$. 
This linear growth is coherent with the behaviour of a finite-time gradient blow-up on a finite grid, where clearly a mesh of spacing $\Delta x$ can only resolve slopes up to order $1/\Delta x$. As a consequence, the largest attainable gradient rises in step as the grid is refined.
As a result $M_{\rm tol}$ must be carefully chosen. In particular, with the value chosen in Test 1 and Test 2 ($M_{\rm tol}=100$) the breakdown is detected at every resolution, whereas with $M_{\rm tol}=1000$ exceeds the $N=400$ ceiling.

Overall, across all three parameters the detected blow-up time remains within the interval $T_{\rm num}\approx0.08$-$0.13$, about one order of magnitude below the theoretical upper bound $T_{\rm UB}$. This provides numerical evidence that the finite-time breakdown of the supercritical regime is a robust numerical fact, and it is not an artifact of the discretization and of the scheme considered.

\subsection{Test 2: effect of chemotaxis on blow-up dynamics}
In Test~2 we investigate numerically the effect of the chemotactic strength on the steepening of the velocity gradient and on the detected breakdown time. We keep fixed the compressive initial datum $A=3.30$ used in Test~1, together with all remaining parameters, and vary
$
\gamma\in\{0,0.5,200,1000\}.
$
We stress that the sufficient analytical threshold depends on $\gamma$. Therefore, the same initial datum need not satisfy the breakdown criterion for all the values of $\gamma$ considered here.

Table~\ref{tab:results_test2} reports the numerical blow-up times obtained for the four different values of the parameter $\gamma$ considered.

\begin{table}[H]
\centering
\caption{Test 2: numerical blow-up times as functions of $\gamma$.}
\label{tab:results_test2}
\begin{tabular}{lcccc}
\toprule
$\gamma$      & $0$      & $0.5$    & $200$    & $1000$   \\
\midrule
$T_{\rm num}$ & $0.0950$ & $0.0950$ & $0.0847$ & $0.0665$ \\
\bottomrule
\end{tabular}
\end{table}

First, we observe that the blow-up time $T_{\rm num}$ corresponding to no chemotaxis affecting the evolution of $u$ ($\gamma=0$), is identical to a small value $\gamma=0.5$. This is somehow expected since, being $-\partial_x u_0(0)=10.39$ far exceeding the threshold ($\mathcal{R}=3.26$ if $\gamma=0$ and $\mathcal{R}=1.54$ if $\gamma=0.5$), the initial compression alone already drives breakdown on a time scale $\mathcal{O}(10^{-2})$, so the moderate chemotaxis does not further accelerate it.

\begin{figure}[H]
\centering
\subfigure{%
\includegraphics[width=0.47\textwidth]{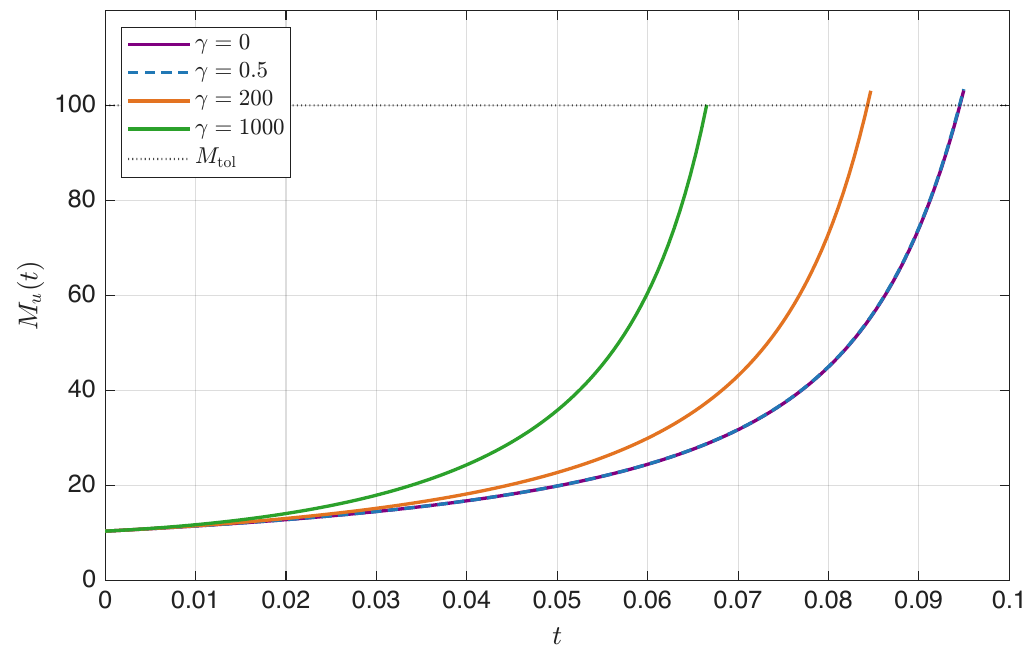}
}
\subfigure{%
\includegraphics[width=0.47\textwidth]{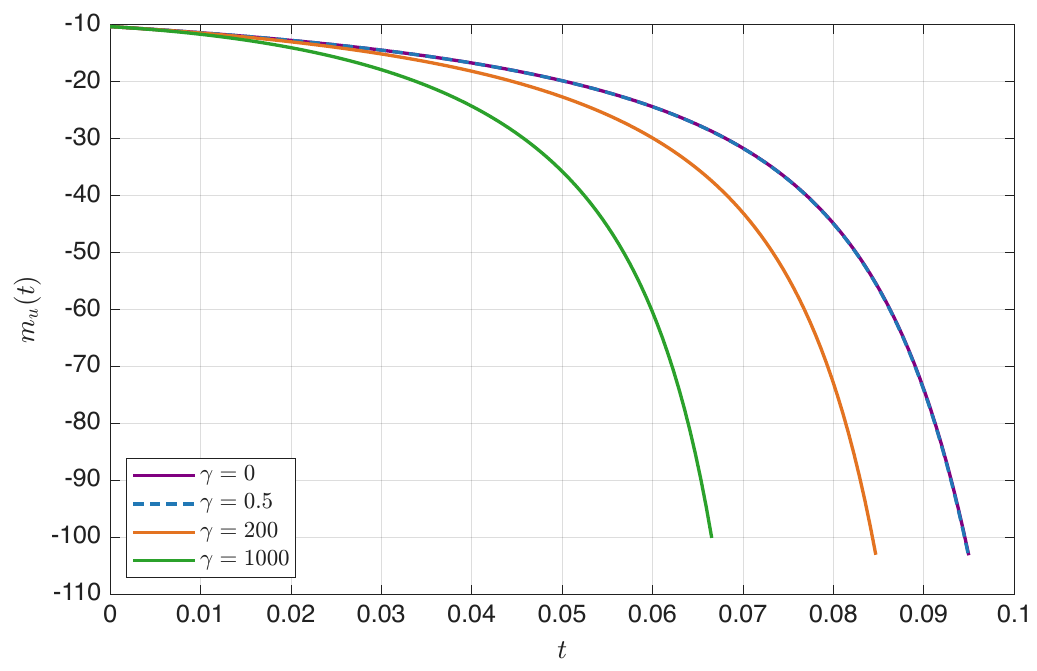}
}
\caption{%
Plot of $m_u(t)$ and $M_u(t)$ for all four values of $\gamma$, for the same compressive initial datum $A_{\rm sup}=3.30$.
Blow-up time decreases with $\gamma$, while the difference between $\gamma=0$ and $\gamma=0.5$ is negligible on the blow-up timescale.}
\label{fig:Mu_minu_seno_variogamma}
\end{figure}

We hence strees the values considered for $\gamma$ in the analysis. For the large values $\gamma=200$ and $\gamma=1000$, $T_{\rm num}$ decreases markedly.
Numerically, $T_{\rm num}$ is decreasing in $\gamma$, suggesting that for strong coupling the magnitude of the chemotactic forcing sets the blow-up time scale.
These results are illustrated in Figure \ref{fig:Mu_minu_seno_variogamma}, which shows the evolution of $M_u(t)$ and $m_u(t)$ up to the blow-up times associated with the different values of $\gamma$.

\begin{figure}[H]
\centering
\subfigure[Density $\rho(x,t)$, $\gamma=200$]{%
  \includegraphics[width=0.47\textwidth]{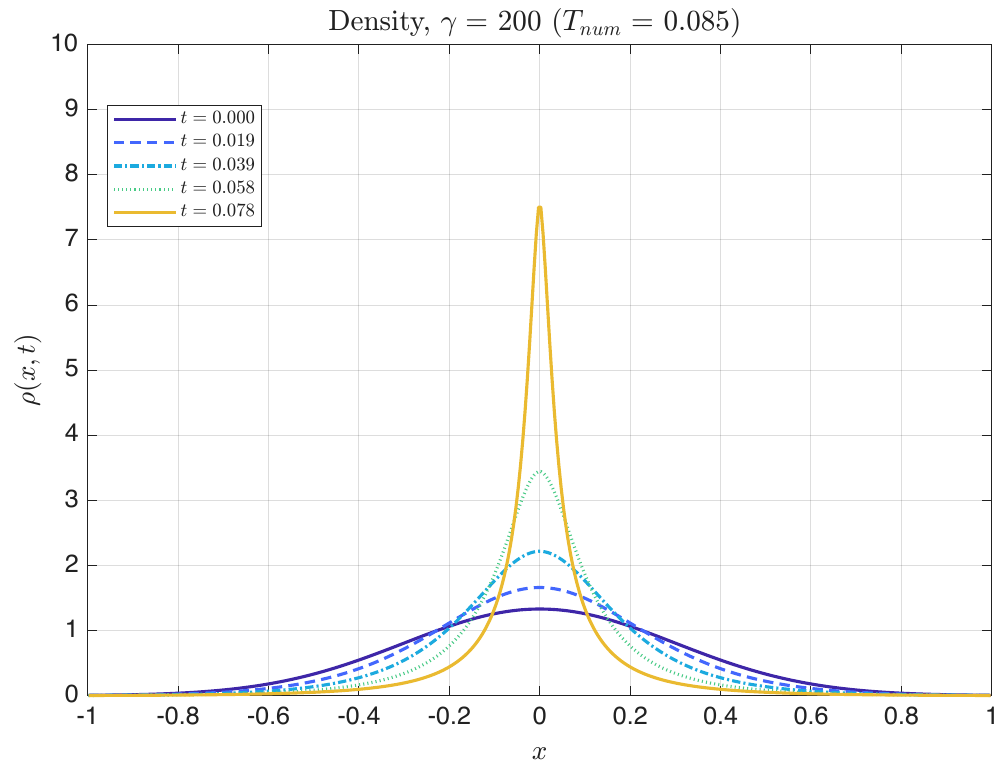}}\hfill
\subfigure[Velocity $u(x,t)$, $\gamma=200$]{%
  \includegraphics[width=0.47\textwidth]{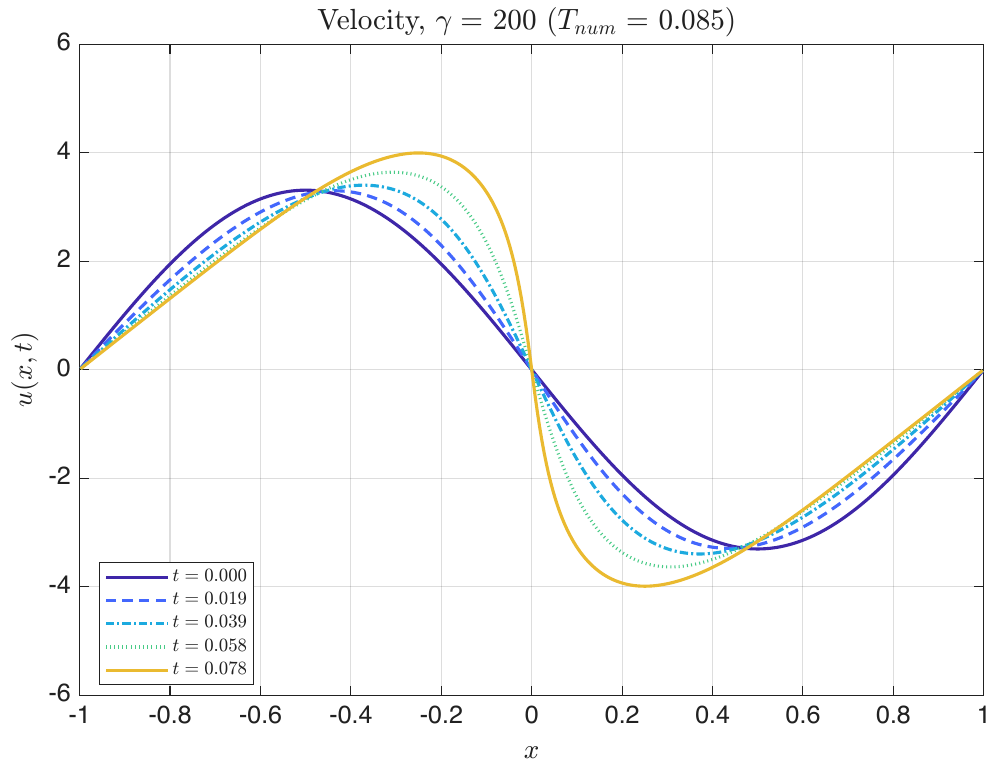}}\\[6pt]
\subfigure[Density $\rho(x,t)$, $\gamma=1000$]{%
  \includegraphics[width=0.47\textwidth]{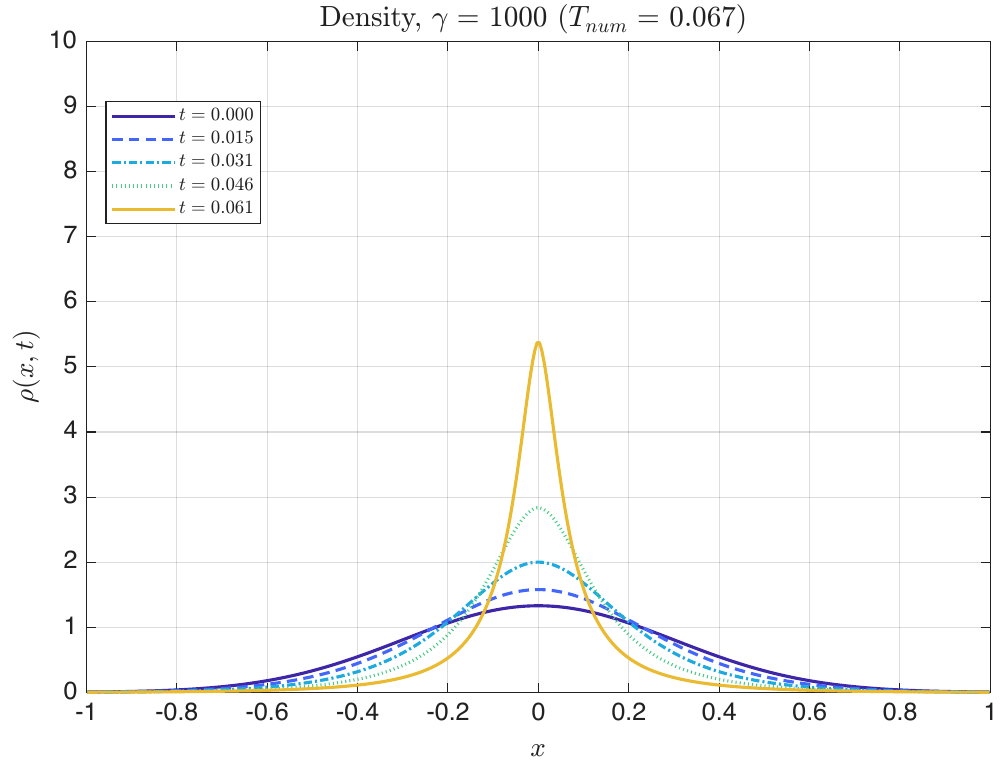}}\hfill
\subfigure[Velocity $u(x,t)$, $\gamma=1000$]{%
  \includegraphics[width=0.47\textwidth]{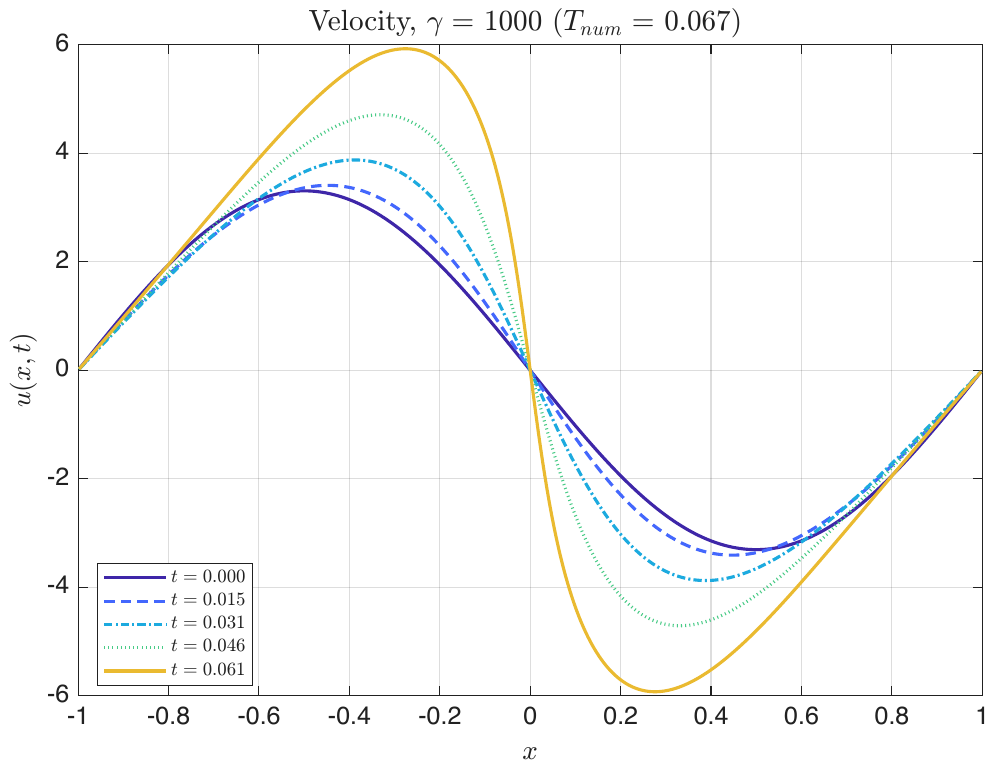}}
\caption{%
Density and velocity profiles evolution.
(a)-(b): $\gamma=200$, $T_{\rm num}=0.0847$; (c)-(d): $\gamma=1000$, $T_{\rm num}=0.0665$.
The chemotactic forcing accelerates concentration of mass and steepening of
the velocity gradient, reducing the blow-up time.
}
\label{fig:panel_g200_g1000}
\end{figure}

Figure \ref{fig:panel_g200_g1000} highlights the qualitative change in blow-up dynamics as the chemotactic coupling $\gamma$ increases.
For $\gamma=200$, the density peak is sharper and the velocity gradient steeper
at blow-up than for $\gamma=0.5$, consistent with the earlier blow-up time $T_{\rm num}=0.0847$.
At $\gamma=1000$ the effect is even more pronounced: blow-up occurs at
$T_{\rm num}=0.0665$, with a markedly more concentrated density profile and
a steeper velocity front near $x=0$.

\section{Conclusion}
\label{sec:conclusion}

In this paper we studied a novel pressureless Euler-alignment system coupled with chemotaxis, generalizing the structure of models derived as the macroscopic counterpart of the class of \textit{hybrid-coupled models}. 
From the analytical viewpoint, under suitable conditions, we first establish local well-posedness and a continuation criterion.
In one dimension, we derive a sufficient condition ensuring continuation beyond a prescribed finite time, while in multiple dimensions, under an additional sign condition along a characteristic, we obtain a sufficient finite-time breakdown criterion.
The core of the work is the study, both analytically and numerically, of blow-up behaviors.
In particular, we derived two complementary criteria.
A sufficient condition for finite-time breakdown in every space dimension and accompanied by an explicit upper bound on the lifespan of the solution is obtained. Moreover, a sufficient condition for global-in-time continuation, established in the one-dimensional setting is provided.
We further clarified the role of chemotaxis through a direct comparison with the pure alignment case, linking our work with previous works of the literature in the field on pressureless alignment systems. Consequently, the subcritical or supercritical feature of the dynamics is no longer determined by the initial configuration alone, but also by the cumulative, memory-type contribution of the chemotactic coupling.
In the constant-coefficient case, the results were also validated numerically. A dedicated sensitivity study shows that the detected blow-up time is robust with respect to mesh size, time step and detection tolerance. Finally, varying the chemotactic strength reveals that a sufficiently large $\gamma$ accelerates the breakdown.
 
Further developments include the treatment of non-constant coefficients, the search for sharp critical thresholds, and the analysis of the long-time behaviour of solutions in the subcritical regime.

\bibliographystyle{plain}
\bibliography{references}

@book{dudley2018real,
  title={Real analysis and probability},
  author={Dudley, Richard M},
  year={2018},
  publisher={Chapman and Hall/CRC}
}

@article{Haskovec2013FlockingDA,
  title={Flocking dynamics and mean-field limit in the Cucker–Smale-type model with topological interactions},
  author={Haskovec, J.},
  journal={Physica D},
  year={2013},
  volume={261},
  pages={42-51},
}

@article {DP19,
    AUTHOR = {Degond, P. and Pulvirenti, M.},
     TITLE = {Propagation of chaos for topological interactions},
   JOURNAL = {Ann. Appl. Probab.},
  FJOURNAL = {The Annals of Applied Probability},
    VOLUME = {29},
      YEAR = {2019},
    NUMBER = {4},
     PAGES = {2594--2612},
      ISSN = {1050-5164,2168-8737},
   MRCLASS = {60K35 (70K45 82B40)},
  MRNUMBER = {3984258},
       DOI = {10.1214/19-AAP1469},
       URL = {https://doi.org/10.1214/19-AAP1469},
}

@article{menci2026numerical,
  title={Numerical Modeling of Flocking Dynamics with Topological Interactions},
  author={Menci, Marta and Paul, Thierry and Rossi, Stefano and Tenna, Tommaso},
  journal={Communications in Computational Physics},
  volume={39},
  number={1},
  pages={240--260},
  year={2026}
}

@book{ambrosio2008gradient,
  author    = {Ambrosio, Luigi and Gigli, Nicola and Savar{\'e}, Giuseppe},
  title     = {Gradient Flows in Metric Spaces and in the Space of Probability Measures},
  edition   = {2nd},
  series    = {Lectures in Mathematics ETH Z{\"u}rich},
  publisher = {Birkh{\"a}user},
  address   = {Basel},
  year      = {2008}
}

@article{albi2018pressureless,
  title={Pressureless Euler alignment system with control},
  author={Albi, Giacomo and Choi, Young-Pil and Haeck, Axel-Stefan},
  journal={Mathematical Models and Methods in Applied Sciences},
  volume={28},
  number={09},
  pages={1635--1664},
  year={2018},
  publisher={World Scientific}
}

@article{menci2026numericalMMS,
  title={Numerical Study on a Multidimensional Pressureless Euler-Type Model with Nonlocal Interactions and Chemotaxis for Collective Cell Migration},
  author={Menci, Marta and Natalini, Roberto and Tenna, Tommaso},
  journal={Multiscale Modeling \& Simulation},
  volume={24},
  number={1},
  pages={267--289},
  year={2026},
  publisher={SIAM}
}

@article{natalini2021mean,
	title={The Mean-Field limit for hybrid models of collective motions with chemotaxis},
	author={Natalini, Roberto and Paul, Thierry},
	journal={SIAM Journal on Mathematical Analysis},
	volume={55},
	number={2},
	pages={900--928},
	year={2023},
	publisher={SIAM}
}

@article{menci2023microscopic,
  title={Microscopic, kinetic and hydrodynamic hybrid models of collective motions with chemotaxis: a numerical study},
  author={Menci, Marta and Natalini, Roberto and Paul, Thierry},
  journal={Mathematics and Mechanics of Complex Systems},
  volume={12},
  number={1},
  pages={47--83},
  year={2024}
}

@article{di2015hybrid,
	title={A hybrid mathematical model for self-organizing cell migration in the zebrafish lateral line},
	author={Di Costanzo, Ezio and Natalini, Roberto and Preziosi, Luigi},
	journal={Journal of Mathematical Biology},
	volume={71},
	number={1},
	pages={171--214},
	year={2015},
	publisher={Springer}
}

@article{costanzo2020hybrid,
	title={A hybrid model of collective motion of discrete particles under alignment and continuum chemotaxis.},
	author={Di Costanzo, Ezio and Menci, Marta and Messina, Eleonora and Natalini, Roberto and Vecchio, Antonia},
	journal={Discrete \& Continuous Dynamical Systems-Series B},
	volume={25},
	number={1},
        pages={443--472},
	year={2020}
}

@article{global2D_Tadmor2017,
  title={Global regularity of two-dimensional flocking hydrodynamics},
  author={He, Siming and Tadmor, Eitan},
  journal={Comptes Rendus. Math{\'e}matique},
  volume={355},
  number={7},
  pages={795--805},
  year={2017}
}

@article{carrillo2016critical,
  title={Critical thresholds in 1D Euler equations with non-local forces},
  author={Carrillo, Jos{\'e} A and Choi, Young-Pil and Tadmor, Eitan and Tan, Changhui},
  journal={Mathematical Models and Methods in Applied Sciences},
  volume={26},
  number={01},
  pages={185--206},
  year={2016},
  publisher={World Scientific}
}

@article{carrillo2017review,
  title={A review on attractive-repulsive hydrodynamics for consensus in collective behavior},
  author={Carrillo, Jos{\'e} A and Choi, Young-Pil and Perez, Sergio P},
  journal={Active Particles, Volume 1: Advances in Theory, Models, and Applications},
  pages={259--298},
  year={2017},
  publisher={Springer}
}

@article{carrillo2026lagrangian,
  title={Lagrangian formulation and Eulerian closure in alignment dynamics},
  author={Carrillo, Jos{\'e} A and Choi, Young-Pil and Tadmor, Eitan},
  journal={arXiv preprint arXiv:2604.10253},
  year={2026}
}

@article{choi2026finite,
  title={Finite-time breakdown of the Euler-alignment system for supercritical initial data},
  author={Choi, Young-Pil and Tadmor, Eitan},
  journal={arXiv preprint arXiv:2607.08060},
  year={2026}
}

@article{tadmor2003critical,
  title={Critical thresholds in 2D restricted Euler-Poisson equations},
  author={Tadmor, Eitan and Liu, Hailiang},
  journal={SIAM Journal on Applied Mathematics},
  volume={63},
  number={6},
  pages={1889--1910},
  year={2003},
  publisher={SIAM}
}

@article{luan2026euler,
  title={On the Euler--Poisson equations with variable background states and nonlocal velocity alignment},
  author={Luan, Kunhui and Tan, Changhui and Wu, Qiyu},
  journal={SIAM Journal on Mathematical Analysis},
  volume={58},
  number={3},
  pages={2251--2290},
  year={2026},
  publisher={SIAM}
}

@article{bhatnagar2020critical,
  title={Critical thresholds in 1D pressureless Euler-Poisson systems with variable background},
  author={Bhatnagar, Manas and Liu, Hailiang},
  journal={Physica D: Nonlinear Phenomena},
  volume={414},
  pages={132728},
  year={2020},
  publisher={Elsevier}
}

@article{tadmor2008global,
  title={On the global regularity of subcritical Euler--Poisson equations with pressure},
  author={Tadmor, Eitan and Wei, Dongming},
  journal={Journal of the European Mathematical Society},
  volume={10},
  number={3},
  pages={757--769},
  year={2008}
}

@article{wei2012critical,
  title={Critical thresholds in multi-dimensional Euler-Poisson equations with radial symmetry},
  author={Wei, Dongming and Tadmor, Eitan and Bae, Hantaek},
  journal={Commun. Math. Sci},
  volume={10},
  number={1},
  pages={75--86},
  year={2012}
}

@article{carrillo2021mean,
  title={Mean-field limits: from particle descriptions to macroscopic equations},
  author={Carrillo, Jos{\'e} A and Choi, Young-Pil},
  journal={Archive for Rational Mechanics and Analysis},
  volume={241},
  number={3},
  pages={1529--1573},
  year={2021},
  publisher={Springer}
}

@incollection{golse2016dynamics,
  title={On the dynamics of large particle systems in the mean field limit},
  author={Golse, Fran{\c{c}}ois},
  booktitle={Macroscopic and large scale phenomena: coarse graining, mean field limits and ergodicity},
  pages={1--144},
  year={2016},
  publisher={Springer}
}

@article{jabin2014review,
  title={A review of the mean field limits for Vlasov equations},
  author={Jabin, Pierre-Emmanuel},
  journal={Kinetic and Related models},
  volume={7},
  number={4},
  pages={661--711},
  year={2014},
  publisher={Kinetic and Related Models}
}

@article{choi2025alignment,
  title={Alignment with Nonlinear Velocity Couplings: Collision Avoidance and Micro-to-Macro Mean-Field Limits},
  author={Choi, Young-Pil and Fabisiak, Micha{\l} and Peszek, Jan},
  journal={SIAM Journal on Mathematical Analysis},
  volume={57},
  number={5},
  pages={5791--5820},
  year={2025},
  publisher={SIAM}
}

@incollection{choi2017emergent,
  title={Emergent dynamics of the Cucker--Smale flocking model and its variants},
  author={Choi, Young-Pil and Ha, Seung-Yeal and Li, Zhuchun},
  booktitle={Active Particles, Volume 1: Advances in Theory, Models, and Applications},
  pages={299--331},
  year={2017},
  publisher={Springer}
}

@book{friedman2008partial,
  title={Partial differential equations of parabolic type},
  author={Friedman, Avner},
  year={2008},
  publisher={Courier Dover Publications}
}

@article{di2018discrete,
	title={{A discrete in continuous mathematical model of cardiac progenitor cells formation and growth as spheroid clusters (Cardiospheres)}},
	author={Di Costanzo, Ezio and Giacomello, Alessandro and Messina, Elisa and Natalini, Roberto and Pontrelli, Giuseppe and Rossi, Fabrizio and Smits, Robert and Twarogowska, Monika},
	journal={Mathematical Medicine and Biology: a Journal of the IMA},
	volume={35},
	number={1},
	pages={121--144},
	year={2018},
	publisher={Oxford University Press}
}

@article{tadmor2014critical,
  title={Critical thresholds in flocking hydrodynamics with non-local alignment},
  author={Tadmor, Eitan and Tan, Changhui},
  journal={Philosophical Transactions of the Royal Society A: Mathematical, Physical and Engineering Sciences},
  volume={372},
  number={2028},
  pages={20130401},
  year={2014},
  publisher={The Royal Society Publishing}
}

@article{menci2019global,
  title={Global solutions for a path-dependent hybrid system of differential equations under parabolic signal},
  author={Menci, Marta and Papi, Marco},
  journal={Nonlinear Analysis},
  volume={184},
  pages={172--192},
  year={2019},
  publisher={Elsevier}
}

@article{menci2022existence,
	title={Existence of solutions for hybrid systems of differential equations under exogenous information with discontinuous source term},
	author={Menci, Marta and Papi, Marco},
	journal={Nonlinear Analysis},
	volume={221},
	pages={112885},
	year={2022},
	publisher={Elsevier}
}

@article{menci2023existence,
  title={Existence of solutions for coupled hybrid systems of differential equations for microscopic dynamics and local concentrations},
  author={Menci, Marta and Papi, Marco and Smarrazzo, Flavia},
  journal={Communications on Pure and Applied Analysis},
  volume={22},
  number={7},
  pages={2146--2168},
  year={2023},
  publisher={Communications on Pure and Applied Analysis}
}

@article{menci2023coupled,
  title={On a coupled hybrid system of nonlinear differential equations with a nonlocal concentration},
  author={Menci, Marta and Papi, Marco and Porzio, Maria Michaela and Smarrazzo, Flavia},
  journal={Journal of Differential Equations},
  volume={361},
  pages={288--338},
  year={2023},
  publisher={Elsevier}
}

\end{document}